\documentclass [12pt,reqno]{article}
\usepackage {graphicx}

\usepackage {amsmath,amssymb,amsfonts,amsthm}

\usepackage {dcolumn}
\newcolumntype {d}[1]{D{.}{.}{#1}}

\usepackage{algorithm}
\usepackage{algorithmicx}
\usepackage{algpseudocode}

\theoremstyle{plain}
\newtheorem {theorem}{Theorem}

\newtheorem {lemma}[theorem]{Lemma}
\newtheorem {proposition}[theorem]{Proposition}
\theoremstyle{definition}

\newtheorem {conjecture}[theorem]{Conjecture}

\usepackage[usenames]{color}
\definecolor{webgreen}{rgb}{0,.5,0}
\definecolor{webbrown}{rgb}{.6,0,0}
\usepackage [colorlinks=true,
linkcolor=webgreen,
filecolor=webbrown,
citecolor=webgreen]{hyperref}

\renewcommand {\epsilon}{\varepsilon}
\renewcommand {\theta}{\vartheta}

\ifxetex
\renewcommand {\C}{\mathbb {C}}
\else
\newcommand {\C}{\mathbb {C}}
\fi
\newcommand {\N}{\mathbb {N}}
\newcommand {\Z}{\mathbb {Z}}
\newcommand {\Q}{\mathbb {Q}}
\newcommand {\Oc}{\mathcal {O}}
\newcommand {\OK}{\Oc_K}
\newcommand {\pf}{\mathfrak {p}}

\newcommand {\Norm}{\operatorname {N}}
\newcommand {\Sl}{\operatorname {Sl}}

\newcommand{\seqnum}[1]{\href{http://oeis.org/#1}{\underline{#1}}}

\begin {document}

\begin{center}
\vskip 1cm{\LARGE\bf Short Addition Sequences for Theta Functions}
\vskip 1cm
\large
Andreas Enge\\
INRIA -- LFANT\\
CNRS -- IMB -- UMR 5251\\
Universit\'e de Bordeaux\\
33400 Talence \\
France \\
\href {mailto:andreas.enge@inria.fr}{\tt andreas.enge@inria.fr} \\
\ \\
William Hart\\
Technische Universit\"at Kaiserslautern\\
Fachbereich Mathematik\\
67653 Kaiserslautern\\
Germany\\
\href {mailto:goodwillhart@googlemail.com}{\tt goodwillhart@googlemail.com}\\
\ \\
Fredrik Johansson\\
INRIA -- LFANT\\
CNRS -- IMB -- UMR 5251\\
Universit\'e de Bordeaux\\
33400 Talence \\
France \\
\href {mailto:fredrik.johansson@inria.fr}{\tt fredrik.johansson@inria.fr}
\end{center}

\vskip .2 in

\begin {abstract}
The main step in numerical evaluation of
classical $\Sl_2(\Z)$ modular forms
is to compute the sum of the first $N$ nonzero terms in the sparse $q$-series
belonging to the Dedekind eta function or the Jacobi theta constants.
We construct short addition sequences to perform
this task using $N + o(N)$ multiplications.
Our constructions rely on the representability of
specific quadratic progressions of integers as sums
of smaller numbers of the same kind. For example, we show that every
generalized pentagonal number $c \geq 5$ can be written as $c = 2a+b$
where $a, b$ are smaller generalized pentagonal numbers.
We also give a baby-step giant-step algorithm
that uses $O(N / \log^r N)$ multiplications for any $r > 0$,
beating the lower bound of $N$ multiplications required
when computing the terms explicitly. These results lead
to speed-ups in practice.
\end {abstract}

\section{Motivation and main results}

We consider the problem of numerically evaluating
a function given by a power series
$f(q) = \sum_{n=0}^{\infty} c_n q^{e_n}$,
where the \emph{exponent sequence}
$E = \{e_n\}_{n=0}^{\infty}$ is a strictly increasing sequence of
natural numbers.
If $|c_n| \leq c$ for all $n$ and the given argument $q$
satisfies $|q| \leq 1 - \delta$
for some fixed $\delta > 0$, then
the truncated series taken over the exponents $e_n \leq T$
gives an approximation of $f(q)$ with error at most
$c |q|^{T + 1} / \delta$, which is accurate to $\Omega(T)$ digits.
This brings us to the question of
how to evaluate the finite sum
$f(q) \approx \sum_{e_n \leq T} c_n q^{e_n}$
as efficiently as possible.
To a first approximation, it is reasonable to attempt to
minimize the total number of multiplications, including coefficient
multiplications and multiplications by $q$.

Our work is motivated by the case $c_n$, $q \in \C$, but
most statements transfer to rings such
as $\Q_p$ and $\C [[t]]$.
The abstract question of how to evaluate truncated power series,
that is, polynomials, with a minimum number of multiplications
may even be asked for an arbitrary coefficient ring. We review a number
of generic approaches in \S\ref {sec:power}.

More precisely, we are interested in highly structured exponent sequences,
namely, sequences given by values of specific quadratic polynomials that
belong to the Jacobi theta constants and to the Dedekind eta function.
Exploiting this structure, one may hope to obtain more efficient algorithms.

The general one-dimensional \emph {theta function} is given by
\begin{equation}
\label {eq:theta}
\theta (\tau, z)
= \sum_{n \in \Z} e^{\pi i n^2 \tau + 2 \pi i n z}
= \sum_{n \in \Z} w^n q^{n^2}
\end{equation}
with $q = e^{\pi i \tau}$ and $w = e^{2 \pi i z}$
for $z \in \C$ and $\tau$ in the upper complex half-plane, that is,
$\Im (\tau) > 0$.
For some $\ell \in \Z_{\geq 1}$ and $a$, $b \in \frac {1}{\ell} \Z$,
the theta function of \emph {level}~$\ell$ and with
\emph {characteristic}~$(a, b)$ is defined as
\[
\theta_{a, b} (\tau, z)
= e^{\pi i a^2 \tau + 2 \pi a (z + b)} \theta (\tau, z + a \tau + b)
= \sum_{n \in \Z} e^{\pi i (n + a)^2\tau + 2 \pi i (n + a)(z + b)}.
\]
The functions of level~$2$ are the classical Jacobi theta functions.
Of special interest are the \emph {theta constants}, the functions of~$\tau$
in which one has fixed $z = 0$ or $w = 1$, respectively, and in particular,
those of level~$2$, given by
\begin {align}
\nonumber
\theta_0 (\tau)
& = \theta_{0, 0} (\tau, 0) = \sum_{n \in \Z} q^{n^2}
= 1 + 2 \sum_{n = 1}^\infty q^{n^2}
= 1 + 2 q \sum_{n = 1}^\infty q^{n^2 - 1}, \\
\label {eq:thetatri}
\theta_1 (\tau)
& = \theta_{0, \frac {1}{2}} (\tau, 0)
= 1 + 2 \sum_{n = 1}^\infty (-1)^n q^{n^2}, \\
\nonumber
\theta_2 (\tau)
& = \theta_{\frac {1}{2}, 0} (\tau, 0)
= 2 \sum_{n = 1}^\infty q^{\frac {1}{4} (2 n + 1)^2}
= 2 q^{\frac {1}{4}} \sum_{n = 1}^\infty q^{n (n + 1)}.
\end {align}
(Here and in the following, when $q$ is defined as $q = e^\gamma$, by a
slight abuse of notation we write $q^{\frac {1}{\ell}}$ for the then
unambiguously defined $e^{\gamma / \ell}$.)
The remaining function $\theta_{\frac {1}{2}, \frac {1}{2}} (\tau, 0)$
is identically~$0$.

Different notational conventions are often used in the literature;
the functions we have denoted by $\theta_0, \theta_1, \theta_2$
are sometimes denoted $\theta_3, \theta_4, \theta_2$ and often
with different factors $\tfrac{1}{2}$ or $\pi$ among the arguments.
Higher-dimensional theta functions are the objects of choice for studying
higher-di\-men\-sional abelian varieties
\cite {Mumford83,Mumford84,Mumford91}.

In dimension~$1$, that is, in the context of elliptic curves,
the \emph {Dedekind eta function} is often more convenient
\cite {Weber08,Schertz02,EnSc04,EnSc13,EnMo14}.
It is a modular form of weight~$\frac {1}{2}$ and level~$24$ defined by
\begin {align}
\nonumber
\eta (\tau)
& = q^{\frac {1}{24}} \prod_{n = 1}^\infty (1 - q^n)
= q^{\frac {1}{24}} \sum_{n=-\infty}^{\infty} (-1)^n q^{n(3n-1)/2}
\quad \cite {Euler83} \\
\label {eq:etapenta}
& = q^{\frac {1}{24}} \left(
\sum_{n = 1}^\infty (-1)^n q^{n (3n-1) / 2}
+ \sum_{n=0}^\infty (-1)^n q^{n (3n+1) / 2} \right)
\end {align}
for $q = e^{2 \pi i \tau}$ (where an additional~$2$ appears in the exponent
compared to the definition of $q$ for theta functions).
It is related to theta functions via
$2 \eta (\tau)^3 = \theta_0 (\tau) \theta_1 (\tau) \theta_2 (\tau)$
\cite [\S34, (10) and (11)]{Weber08}, and
$\eta (\tau) = \zeta_{12} \, \theta_{- \frac {1}{6}, \frac {1}{2}} (3 \tau, 0)$
with $\zeta_{12} = e^{2 \pi i / 12}$.
The latter property can be proved easily as an equality of formal series.

Other functions that can be expressed in terms of theta functions
include Eisenstein series $G_{2k}(\tau)$
and the Weierstrass elliptic function $\wp(z; \tau)$.
Theta functions are also useful in physics for solving the heat equation.

Another motivation for looking at theta and eta functions comes from
complex multiplication of elliptic curves. The moduli space of
complex elliptic curves is parameterized by the j-invariant,
given by a $q$-series
$j(\tau) = \frac {1}{q} + 744 + 196884 q + 21493760 q^2 + \cdots$
for $q = e^{2\pi i \tau}$, which can be obtained explicitly from the series
of the theta or eta functions as \cite[\S54, (5); \S34, (11)]{Weber08}
\[
j = \left( \frac {f_1^{24} + 16}{f_1^8} \right)^3
\text { with }
f_1 (\tau) = \frac {\eta (\tau / 2)}{\eta (\tau)},
\]
or \cite[\S34, (10) and (11); \S54, (6); \S21, (14)]{Weber08}
\begin {equation}
j = 32 \, \frac {(\theta_0^8 + \theta_1^8 + \theta_2^8)^3}
                {(\theta_0 \theta_1 \theta_2)^8}.
\end {equation}

The series for $j$ is dense, and the coefficients in front of~$q^n$
asymptotically grow as
$2^{1/2} n^{-3/4} e^{4\pi\sqrt{n}}$ \cite{Rademacher38};
so it is in fact preferable to obtain its values from values of theta
and eta functions: Their sparse series imply that $O (\sqrt T)$ terms
are sufficient for a precision of $\Omega (T)$ digits, and they
furthermore have coefficients $\pm 1$.

Evaluating $j$ at high precision is a building block for the complex
analytic method to
compute ring class fields of imaginary-quadratic number fields and then
elliptic curves with a given endomorphism ring \cite {Enge09cla}, or
modular polynomials encoding isogenies between elliptic curves \cite {Enge09mod}.
For example, the Hilbert class polynomial for a quadratic
discriminant $D < 0$ is given by
\[
H_D(x) =
\prod_{(a,b,c)} \left( x - j\left(\frac{-b+\sqrt{D}}{2a} \right) \right)
\in \Z [x],
\]
where $(a,b,c)$ is taken over the primitive reduced binary
quadratic forms $ax^2 + bxy + cy^2$ with $b^2 - 4ac = D$.
The exact coefficients of $H_D$
can be recovered from $|D|^{1/2+o(1)}$-bit numerical approximations.

The bit complexity of evaluating the theta or eta functions at a precision
of $T$ digits via their $q$-series is in $O (T^{3/2+o(1)})$.
Asymptotically for $T$ tending to infinity there is a quasi-linear
algorithm with bit complexity $O (T \log^2 T \log \log T)$
\cite {Dupont11}; it uses the arithmetic-geometric mean (AGM) iteration
together with Newton iterations on an approximation computed
at low precision by evaluating the series.
The crossover point where the
asymptotically faster quasi-linear algorithm wins is quite high.
In earlier work~\cite [Table~1]{Enge09cla}, it was seen to occur
at a precision of about $250\,000$~bits, used to compute
a class polynomial of size about $5$~GB. So in most practical situations,
series evaluation is faster. This is also due to the experimental
observation, implemented in the software~CM~\cite {cm}, that there
are particularly short addition sequences for the exponents
in the $q$-series of the Dedekind eta function,
which lead to a small constant in the complexity in $O$-notation.

Looking at eta and theta functions, respectively, in \S\S\ref {sec:eta}
and~\ref {sec:theta}, we show that this is not a coincidence, but a
consequence of their structured exponents.

Some of our results depend on the Bateman-Horn conjecture for the special
case of only one polynomial, which can be summarized as follows:
\begin {conjecture}[Bateman-Horn \cite {BaHo62}]
\label {conj:bh}
Let $f \in \Z [X]$ be a polynomial with positive leading coefficient
such that for every prime~$p$, there is an $x$ modulo $p$ with
$p \nmid f (x)$. Then there exists a constant
$C > 0$ such that the number of primes among the first $N$ values
$f (1), f (2), \ldots, f (N)$ is asymptotically equivalent to
$C N / \log N$ for $N \to \infty$.
\end {conjecture}
In other words, the density of primes among the values of $f$ is the same
as the density of primes among all integers of the same size, up to a
correction factor $C$, which is given by an Euler product encoding the
behavior of $f$ modulo primes. The hypothesis of the conjecture
is clearly necessary; if it is not satisfied, then all values of $f$ are
divisible by the same prime~$p$, so the only prime potentially occurring is
$p$ itself, and this can happen only a finite number of times (and then
indeed one of the Euler factors defining~$C$ vanishes).
All polynomials $f$ that we consider have $f (0) = 1$, so the hypothesis
is trivially verified.

In particular, we show the following:
\begin {theorem}
\label {th:addseq}
\quad
\begin {enumerate}
\item
The first $N$ terms of the $\eta$ series may be evaluated with
$N + O (1)$ squarings and $N + O (1)$ additional multiplications.
(This follows from Theorem~\ref {th:etatwoaplusb}.)
\item
The first $N$ terms of a series yielding $\eta$
may be evaluated with $N + O (N / \log N)$ multiplications,
assuming Conjecture~\ref {conj:bh} for the polynomials
$18 n^2 - 6 n + 1$ and $18 n^2 + 6 n + 1$.
(This follows from Theorems~\ref {th:etaaplusb} and~\ref {th:etatwoaplusb}.)

The same assertion holds for $\theta_0$ and $\theta_1$
assuming Conjecture~\ref {conj:bh} for the polynomials
$4 n^2 + 1$ and $2 n^2 + 2 n + 1$.
(This follows from Theorems~\ref {th:thetazero} and~\ref {th:thetazerothree}.)

The same assertion holds for $\theta_2$
assuming Conjecture~\ref {conj:bh} for the polynomial
$2 n^2 + 2 n + 1$.
(This follows from Theorems~\ref {th:thetatwo} and~\ref {th:thetatwothree}.)
\item
Truncating $\theta_0$, $\theta_1$ and $\theta_2$ to $N$ terms
each, only $2 N$ monomials occur.
The first $N$ terms of series yielding
all three theta constants (in the same argument $q$) may
be evaluated with $2 N + O (1)$ multiplications.
(This is Theorem~\ref {th:thetaaplusb}.)
\end {enumerate}
\end {theorem}

The number of multiplications needed to evaluate a series is closely related
to the number of additions needed to compute the values of its exponent sequences,
as we discuss further in \S\ref {sec:power}.
Dobkin and Lipton have previously proved
a lower bound of $N + N^{2/3 - \epsilon}$ additions
for computing the values of certain polynomials at the first
$N$ integers~\cite{DoLi80}, which in particular holds for the squares occurring as exponents
of~$\theta_0$. Dobkin and Lipton conjecture that this lower
bound holds for
arbitrary (non-linear) polynomials. While not exactly a counterexample,
the third point of Theorem~\ref {th:addseq} shows that the conjecture
does not hold when the values of two polynomial sequences are
interleaved.

Finally in \S\ref {sec:bsgs} we present a new baby-step giant-step algorithm
for evaluating theta or eta functions that is asymptotically faster than
any approach computing all monomials occurring in the truncated series.
\begin {theorem}
\label {th:bsgs}
There is an effective constant $c > 0$ such that the series for
$\eta$, $\theta_0$, $\theta_1$ or $\theta_2$, truncated at $N$ terms,
is evaluated by the baby-step giant-step algorithm of \S\ref {sec:bsgs}
with less than
$N^{1 - c / \log \log N}$ multiplications.
(This is a consequence of Theorem~\ref {th:quadpolym}.)
\end {theorem}

Though asymptotically not as fast as the AGM method, this
algorithm gives a speed-up in the practically
interesting range from around $10^3$ to $10^6$ bits,
and further raises the crossover point for the AGM method;
see \S\ref {sec:bench}.

The baby-step giant-step algorithm relies on finding a suitable sequence
of parameters~$m$ such that the exponent sequence
takes few distinct values modulo~$m$; we solve this problem
for general quadratic polynomials and explicitly describe the
parameters $m$ corresponding to the squares, trigonal and pentagonal
numbers occurring as exponents of the eta and theta functions.

The general theta series~\eqref{eq:theta} can be seen as the
Laurent series
$\sum_{n \in \Z} f_n w^n$. 
Theorem~\ref{th:addseq} implies a fast way to compute the coefficients $f_n$.
This speeds up computing the theta function~\eqref{eq:theta} for general~$q, w$
and consequently also speeds up computing
elliptic functions and $\Sl_2 (\Z)$ modular forms via theta functions.
The baby-step giant-step algorithm of
Theorem~\ref{th:bsgs} does not compute the coefficients $f_n$ explicitly.
It speeds up modular forms further, but this speed-up only
applies to the special case $w = 1$ (or other simple algebraic values of $w$,
by a slight generalization), so it is less useful for elliptic functions.

\section {Power series and addition sequences}
\label {sec:power}

In this section, we review some of the known techniques for evaluating
a truncated power series
$f(q) = \sum_{n=0}^{N} c_n q^{e_n}$, where the exponent sequence
$(e_n)_{n = 0}^\infty$ is strictly increasing and the cut-off parameter
$N$ is chosen such that $e_N \leq T$ and $e_{N+1} > T$ for some truncation order $T$
depending on the required precision. (As mentioned before, if a lower
bound on $|q|$ is given, $T$ will be linear in the desired bit precision.)
We let $E = (e_n)_{n=1}^N$ and distinguish between the cases where this sequence
is dense or sparse.

\subsection{Dense exponent sequences}

If the exponent sequence $E$ is dense, that is, $N \in \Omega (T)$,
then Horner's rule is optimal in general.
For example, if $E$ is an arithmetic progression with
step length $r$, then $T/r + O(\log r)$ multiplications suffice.

It is possible to do better if the coefficients $c_n$ have a special form.
Of particular interest is when a multiplication $c_n \cdot q$ is
``cheap'' while a multiplication such as $q \cdot q$ is ``expensive''.
This is the case, for instance, when $c_n$ are small integers or rationals with
small numerators and denominators and~$q$ is a high-precision
floating-point number.
In this case we refer to $c_n$ ``scalars'' and call $c_n \cdot q$
a ``scalar multiplication'' while $q \cdot q$ is called a ``nonscalar multiplication''.
All $N$ multiplications in Horner's rule are nonscalar.
Paterson and Stockmeyer introduced a baby-step giant-step algorithm
that reduces the number of nonscalar multiplications~\cite{PaSt73} to $\Theta(N^{1/2})$,
originally
for the purpose of evaluating polynomials of a matrix argument (which explains
the ``scalar'' terminology).
The idea is to write the series as
\begin{equation}
\sum_{n=0}^{N-1} c_n q^n
= \sum_{k=0}^{\lceil N/m \rceil - 1} \left(q^m\right)^k
\left( \sum_{j=0}^{m-1} c_{mk+j} q^j \right)
\label {eq:bsgsformula}
\end{equation}
for some splitting parameter $m$.
The ``baby-steps'' compute the powers $q^2, q^3, \ldots, q^m$ once and
for all, so that all inner sums may be obtained using multiplications
by scalars. The outer polynomial evaluation with respect to $q^m$ is
then done by Horner's rule using ``giant-steps''.
This requires about $N$ multiplications by scalars and, by choosing
$m \in \Theta (N^{1/2})$ and thus balancing the baby- and giant-steps,
$\Theta (N^{1/2})$ nonscalar multiplications.

For even more special coefficients $c_n$, further techniques exist~\cite{ChCh88,Bernstein2008}:
\begin{itemize}
\item If $E$ is an arithmetic progression and the coefficients $c_n$ satisfy
a linear recurrence relation with polynomial coefficients, then
$N^{1/2+o(1)}$ arithmetic operations (or $N^{3/2+o(1)}$ bit operations)
suffice if fast multipoint evaluation of polynomials is used.
An improved version of the Paterson-Stockmeyer algorithm also exists
for such sequences~\cite{Smith1989,Johansson14}.
\item If both $q$ and the coefficients $c_n$ are scalars of a suitable
type, binary splitting should be used.
For example, if the ``scalars'' are rational numbers (or elements of a fixed number field)
with $O(\log n)$ bits, the bit complexity is reduced to the quasi-optimal $N^{1+o(1)}$.
This result also holds if $E$ is an arithmetic progression,
$q \in \overline \Q$, and the $c_n$ satisfy a linear recurrence relation
with coefficients in $\overline \Q(n)$.
\end{itemize}

The last technique is useful for computing many mathematical
functions and constants, especially those represented by
hypergeometric series, where $q$ often will be algebraic.
It appears to be less useful in connection with theta series,
where $q$ usually will be transcendental.

\subsection{Sparse exponent sequences and addition sequences}

If the exponent sequence $E$ is sparse, for instance if
$e_n \in \Theta(n^{\alpha})$ so that $N \in \Theta (T^{1/\alpha})$
for some $\alpha > 1$, methods designed for dense series
may become inferior to even naively computing separately the powers
of $q$ that are actually needed and evaluating $\sum_n c_n q^{e_n}$ as
written. Addition sequences \cite [Definition~9.32]{CoFrAvDoLaNgVe06}
provide a means of saving work by
computing the needed powers of~$q$ simultaneously.

An \emph{addition sequence} consists of a set of positive integers $A$
containing~$1$,
and for every $c \in A_{>1}$, a pair $a, b \in A_{<c}$
such that $c = a + b$.
An addition sequence $A \supseteq E$ allows us to compute $\{q^e : e \in E\}$
using at most $|A|-1$ multiplications
\[
q^c = q^a \cdot q^b, \quad c \in A.
\]
Given a list of positive integers $E = \{e_1, e_2, \ldots, e_N\}$
with $e_1 < e_2 < \cdots < e_N$, we may have to insert extra elements
to obtain an addition sequence.
For example, the Fibonacci sequence $\{ 1, 2, 3, 5, 8, 13, \ldots \}$
trivially forms an addition sequence without adding more elements, while
the sequence of squares $\{ 1, 4, 9, 16, 25, 36, \ldots \}$ requires
adding intermediate steps.
Minimizing the number of insertions required to
form an addition sequence
becomes an interesting problem; its associated decision
problem is NP-complete in general \cite [Theorem~3.1]{DoLeSe81}.

\begin{algorithm}[hbt]
\caption{Short addition sequence}
\label {alg:addseq}
\begin{algorithmic}
\Require A finite list of positive integers $E$
\Ensure An addition sequence $A \supseteq E$
\State Let $A = E$.
\State While some element $c \in A$, $c \ne 1$,
is not a sum of two smaller elements of~$A$,
insert $\lfloor c / 2 \rfloor$ and $\lceil c / 2 \rceil$ into $A$.
\end{algorithmic}
\end{algorithm}

A straightforward approach, Algorithm~\ref{alg:addseq}, is a close
relative of the double-and-add algorithm for the case of a single exponent,
and it is easy to show that it produces
an addition sequence of length at most $O(N \log e_N) = O (N \log T)$.
In practice, it is observed to produce nearly optimal addition sequences
for reasonably dense input.
A more elaborate method (Yao 1976,
cited by Knuth~\cite[\S4.6.3, exercise 37]{Knuth98}) gives
the upper bound
\[
O\left( N \frac{\log e_N}{\log \log e_N}
+ \log e_N
+ \frac{\log e_N \log \log \log e_N}{(\log \log e_N)^2} \right).
\]

We can improve the upper bounds for sequences of a special form.
For any integer-valued polynomial $f \in \Q [X]$ of degree $d$, the consecutive
values $f(1), f(2), \ldots$ can be computed using $d$ additions
for each new term by the approach of finite differences, letting $f_d = f$
and considering the system of coupled recurrence equations
$f_k(X+1) = f_k(X) + f_{k-1}(X)$, $1 \leq k \leq d$,
in which $\deg(f_k) = k$.

For the quadratic exponent sequences $E$ appearing
in the Dedekind eta function and the Jacobi theta functions,
this implies a cost of two multiplications to generate
each new power $q^{f(n)}$.
We call these the \emph{classical addition sequences},
cf.\ Table~\ref {tab:classical}.

\begin{table}[hbt]
\begin {center}
\begin{tabular}{c|c|c}
$f_2(n)$ & $f_1(n) = f_2(n+1) - f_2(n)$ & $f_0(n) = f_1(n+1) - f_1(n)$ \\
\hline
$n^2$          & $2n+1$ & 2 \\
$n(n+1)$       & $2n+2$ & 2 \\
$n (3n-1) / 2$ & $3n+1$ & 3 \\
$n (3n+1) / 2$ & $3n+2$ & 3 \\
\end{tabular}
\end {center}
\caption {\label {tab:classical} Construction of the
classical addition sequences for squares, trigonal numbers,
and pentagonal numbers via finite differences.}
\end{table}

The classical addition sequences
are often used in implementations~\cite [Algorithm~6.32]{Cohen00},
but they are still not optimal.
For the sequence of squares, Dobkin and Lipton~\cite{DoLi80} give an algorithm
which requires $N + O(N / \sqrt{\log N})$ additions.
Asymptotically, this amounts to a cost of only $1 + o (1)$ multiplications
for each power $q^{n^2}$ in the series for $\theta_0$ or $\theta_1$.
The second point of Theorem~\ref {th:addseq} (heuristically) improves
this bound to $N + O (N / \log N)$.

\subsection{Cost of an addition sequence}

Since squaring is usually cheaper than a general multiplication, it makes
sense to count the number of doublings $c = 2a$ separately from
general additions $c = a+b$ in an addition sequence. We may
even go further and regroup entries in an addition sequence, thus obtaining
more complex atomic operations, to each of which a different cost can be
assigned.

Suppose in particular that multiplying two real floating point numbers
costs~$M$, that squaring such a number costs $S \leq M$ and that additions
and subtractions and, by extension, multiplications by small
integer constants are essentially free. (In fact, we will not need to
consider integer constants other than $1$ and $-1$.)
At high precision, multiplication may
rely on the fast Fourier transform (FFT), the dominant steps of which are the
computation of two forward and one inverse transforms.
When squaring, one of the forward transforms can be skipped, resulting
asymptotically in $S = \frac {2}{3} M$. Using school book multiplication,
one would have $S = \frac {1}{2} M$ asymptotically instead. (We can lower
costs some more by saving the Fourier transform of an operand that is reused
several times, but this results in a more complicated analysis, which we do
not pursue here.)

For complex numbers represented by two reals in Cartesian coordinates,
we have
\begin {align*}
(x+yi)^2     & = (x^2 - y^2) + 2 x \cdot y i, \\
(x+yi)(t+ui) & = (x \cdot t-y \cdot u)
                 + \big( (x+y) \cdot (t+u) - xt - yu \big) i, \\
(x+yi)^3     & =  x\cdot (x^2 - 3y^2) + y\cdot (3x^2 - y^2) i.
\end {align*}

Accordingly, if the complex numbers $q^a$ and $q^b$ have already been
computed (i.e.,\ if $a$ and $b$ are already in the addition sequence), then
we may evaluate the cost for forming the respective new power
(i.e.,\ extending the addition sequence, possibly twice), in increasing
order as in Table~\ref {tab:cost}.

\begin{table}[hbt]
\begin {center}
\begin{tabular}{l|r|l|l}
Step in addition sequence & Generic cost & FFT & School book \\
\hline
$2a$ & $2S + M$ & $2.33 M$ & $2 M$ \\
$a+b$ & $3M$ & $3 M$ & $3 M$ \\
$3a$ & $2S + 2M$ & $3.33 M$ & $3 M$ \\
$4a$ & $4S + 2M$ & $4.67 M$ & $4 M$ \\
$2a+b$ or $2(a+b)$ & $2S + 4M$ & $5.33 M$ & $5 M$
\end{tabular}
\end {center}
\caption {\label {tab:cost}
Costs associated with evaluating complex series using addition sequences.}
\end{table}

\section {Addition sequences for the Dedekind eta function}
\label {sec:eta}

The exponents $e_n = n (3 n -1) / 2$ in \eqref {eq:etapenta} for
$n \geq 1$ are called (ordinary) \emph {pentagonal numbers};
for arbitrary~$n$, \emph {generalized pentagonal numbers} (\seqnum{A001318}
in the On-Line Encyclopedia of Integer Sequences).
In the ordered sequence of exponents, ordinary and
generalized pentagonal numbers alternate.

The sequence of generalized pentagonal numbers is too sparse to
be an addition sequence.
The classical addition sequence effectively doubles the density.
Our observation is that an addition sequence
can be formed by occasionally inserting
an extra doubling (that is, performing an extra squaring
when evaluating the series).

\begin{algorithm}[hbt]
\caption{Dedekind eta function using optimized addition sequence}
\label {alg:eta}
\begin{algorithmic}
\Require $T \ge 2$, $q \in \C$
\Ensure $S = \sum_{e_n \leq T} s_n q^{e_n}$, where
  $E = (e_n)_{n=1}^\infty = (0, 1, 2, 5, 7, \ldots)$
  is the ordered sequence of generalized pentagonal numbers, and
  $s_n \in \{ \pm 1 \}$ is such that $S$ approximates the value of $\eta$
\State $N \gets$ the maximal $n$ such that $e_n \leq T$
\State $S \gets 1 - q, A \gets \{ 1 \}, Q \gets \{ q \}$
\For {$c \gets e_3, \ldots, e_N$}
    \If{$c = 2a$ for some $a \in A$}
        \State $q' \gets (q^a)^2$
    \ElsIf{$c = a+b$ for some $a$, $b \in A$}
        \State $q' \gets (q^a) \cdot (q^b)$
    \ElsIf{$c = 2a+b$ for some $a$, $b \in A$}
        \State $q' \gets (q^a)^2 \cdot q^b$
    \EndIf
    \If {$s_n = +1$}
        \State $S \gets S + q'$
    \Else
        \State $S \gets S - q'$
    \EndIf
    \State $A \gets A \cup \{ c \}$,
           $Q \gets Q \cup \{ q' \}$
\EndFor
\end{algorithmic}
\end{algorithm}

Algorithm~\ref{alg:eta} attempts to write each occurring power of $q$
as a product of previously computed powers. It first attempts the
cheapest operation (squaring)
according to Table~\ref{tab:cost} and proceeds to more expensive operations
if this fails.

In the following, we will prove that the algorithm is correct;
that is, at least one of the branches can always be entered.
In fact, the $c = 2a+b$ case alone is guaranteed to succeed.
That is, every generalized pentagonal number is a sum of
a smaller generalized pentagonal number and twice a smaller
generalized pentagonal number (Theorem~\ref{th:etatwoaplusb}).
We also show that the $c = a+b$ case heuristically almost always succeeds
(Theorem~\ref{th:etaaplusb}),
so that Algorithm~\ref{alg:eta} approaches on average
one multiplication per computed term.

Experimentally, we observe that Algorithm~\ref{alg:eta} uses
slightly fewer multiplications than an addition sequence
constructed with Algorithm~\ref{alg:addseq} when $N$ is large, and a
larger proportion of the multiplications are squarings
(see Figure~\ref{fig:etaefficiency} in \S\ref {ssec:modval}).

The starting point for our considerations is the following
characterization of the generalized pentagonal numbers,
which is immediate from their definition.

\begin {lemma}
\label {lm:penta}
When restricted to generalized pentagonal numbers, the strictly increasing
map
\begin{equation}
\label {eq:sigmadef}
\sigma : c \mapsto \sqrt {24 c + 1}
\end{equation}
is a bijection between generalized pentagonal numbers and positive integers
coprime to~$6$. More precisely, it sends ordinary pentagonal numbers to
integers that are $5 \pmod 6$ and generalized, non-ordinary pentagonal numbers
to integers that are $1 \pmod 6$.
\end {lemma}

\begin {proof}
The equation $c = (3n-1) n / 2$ is equivalent to
$24 c + 1 = (6 n - 1)^2 = (6 (-n) + 1)^2$.
\end {proof}

The first few generalized pentagonal numbers and associated values
of~$\sigma$ are given in Table~\ref {tab:gpn}.

\begin {table}[hbt]
\begin {center}
\begin {tabular}{l|rrrrrrrrrrrrrrr}
$c$          & $0$ & $1$ & $2$ &  $5$ &  $7$ & $12$ & $15$ & $22$ & $26$
             & $35$ & $40$ & $51$ & $57$ & $70$ \\
\hline
$\sigma (c)$ & $1$ & $5$ & $7$ & $11$ & $13$ & $17$ & $19$ & $23$ & $25$
             & $29$ & $31$ & $35$ & $37$ & $41$
\end {tabular}
\end {center}
\caption {\label {tab:gpn}
Generalized pentagonal numbers}
\end {table}

\subsection {One squaring and one multiplication}

The following result provides a proof of the first point of
Theorem~\ref {th:addseq}.

\begin {theorem}
\label {th:etatwoaplusb}
Every generalized pentagonal number $c \geq 5$ is the sum of a smaller one
and twice a smaller one, that is, there are generalized pentagonal numbers
$a, b < c$ such that $c = 2 a + b$.
\end {theorem}

In other words, the series of~$\eta$ may be computed with one multiplication
and one square instead of two multiplications per term, reducing the cost
in the FFT model from $6 M$ to $5.33 M$ according to Table~\ref {tab:cost}.

Hirschhorn shows \cite [(1.20)]{Hirschhorn09}
that the number of ways in which an arbitrary number~$c$ can be written as
twice a generalized pentagonal number plus another pentagonal number is given
by
\[
d_{1,8} (24 c + 3) - d_{7,8} (24 c + 3)
- \big( d_{1,8} ((8 c + 1) / 3) - d_{7,8} ((8 c + 1) / 3) \big),
\]
where $d_{i, j}$ counts the number of positive divisors that are $i\pmod j$
for integral arguments, and equals~$0$ for non-integral rational arguments.
Using quadratic reciprocity and Proposition~\ref {prop:powersofthree} below,
one can show that this quantity is at least~$1$ if $c$ is a generalized
pentagonal number. We prefer to give direct proofs of
Theorem~\ref {th:etatwoaplusb} as well as for similar results below, as they
are instructive and are scarcely more involved than proofs relying on
Hirschhorn's results.

Using Lemma~\ref {lm:penta}, the theorem becomes essentially a statement
about representability of integers as sums of squares. Its proof relies
on the following well-known lemma, for which we give a quick proof
for the sake of self-containedness.

We say that a quadratic form $q (X, Y) = A X^2 + B X Y + C Y^2$
\emph {represents} an integer~$k$ if there are $x$, $y \in \Z$
such that $k = q (x, y)$. The representation is \emph {primitive}
if $x$ and $y$ are coprime.
We are only concerned with the case $B = 0$, and then we say that
the representation is \emph {positive} if $x$, $y > 0$.
If moreover $A = C = 1$, we say that the representation is
\emph {ordered} if $0 < x < y$.

\begin {lemma}
\label {lm:eight}
A positive integer~$k$ is primitively represented by the quadratic
form $2 X^2 + Y^2$ if and only if all its odd prime divisors are~$1$
or $3 \pmod 8$ and it is not divisible by~$4$. Its number of positive
primitive representations is then given by $2^{\omega' (k) - 1}$,
where $\omega' (k)$ denotes the number of odd primes dividing~$k$.
\end {lemma}

Euler \cite{Euler61} proves that a number is primitively representable
in this way if and only if all its prime factors are, and that being
congruent to~$1$ or $3 \pmod 8$ is a necessary condition for odd primes.
Conversely, Euler \cite[p.~628]{Euler55} shows that any prime number that
is congruent to~$1 \pmod 8$ is represented this way.
The missing case of primes congruent to~$3 \pmod 8$ is treated by
Dickson \cite[p.~9]{Dickson23} with a proof attributed to Pierre-Simon
de la Place; we were, however, unable to locate the original reference
of 29 pages from 1776,
\textit {Th\'eorie abr\'eg\'ee des nombres premiers}, in the Gauthier--Villars
edition of the \textit {{\OE}uvres compl\`etes de Laplace} printed in Paris
between 1878 and 1912; the previous and less complete edition of the
\textit {{\OE}uvres de Laplace} printed by the Imprimerie Royale in Paris
between 1843 and 1847 does not contain any number theoretic articles.
Concerning the number of representations, Euler \cite{Euler61} shows that
only the odd primes need to be taken into account, and essentially contains
the formula for square-free~$k$, mentioning explicitly products of two
or three primes and hinting at products of four primes.
Using modern number theory concepts, it is easy to provide a complete
proof of the statements.

\begin {proof}
Representations of~$k$ by $2 X^2 + Y^2$ correspond to elements
$\alpha = x \sqrt {-2} + y$ of the ring of integers
$\OK = \Z [\sqrt {-2}]$ of $K = \Q (\sqrt {-2})$ such that
$\Norm_{K / \Q} (\alpha) = k$.
They are primitive if and only if $\alpha$ is primitive in the sense that
it is not divisible in~$\OK$ by a positive rational integer other than~$1$.
Let $k = 2^{e_0} \prod_{i = 1}^{\omega' (k)} p_i^{e_i}$ be the prime
factorization of~$k$.
A necessary condition for the existence of a primitive representation,
assumed to hold in the further discussion, is that all the $p_i$ are
split in~$K$, which is indeed equivalent to $p_i \equiv 1$ or $3 \pmod 8$~\cite[p.~1]{Cox89},
and that $e_0 \in \{ 0, 1 \}$.
Write $p_i \OK = \pf_i \overline \pf_i$, where $\overline \cdot$ denotes
complex conjugation, the non-trivial Galois automorphism of~$K / \Q$,
with $p_i = \Norm_{K / \Q} (\pf_i)$; and write $2 \OK = \pf_0^2$.
Then $\alpha \in \OK$ is of norm~$k$
(and thus leads to a representation of~$k$) if and only if there are
$\alpha_i \in \{ 0, \ldots, e_i \}$ such that
$\pf_0^{e_0}
\prod_{i=1}^{\omega' (k)} \pf_i^{\alpha_i} \overline \pf_i^{e_i - \alpha_i}$
is a principal ideal generated by~$\alpha$, and the representation is
primitive if and only if none of the $\pf_i$ and $\overline \pf_i$ appear
simultaneously, that is, $\alpha_i \in \{ 0, e_i \}$.
Here the ring $\OK$ is principal, so that principality does not form a
restriction. Letting $\pf_i = \pi_i \OK$ with $\pi_i \in \OK$, the primitive
elements of norm~$k$ are exactly the
\[
\alpha = \epsilon \pi_0^{e_0} \prod_{i=1}^{\omega' (k)} \omega_i^{e_i},
\]
where $\epsilon \in \{ \pm 1 \}$ is a unit in~$\OK$ and
$\omega_i \in \{ \pi_i, \overline \pi_i \}$,
so there are $2^{\omega' (k) + 1}$ of them.
Now there are four possibilities for the signs of $x$ and $y$, meaning
that there are $2^{\omega' (k) - 1}$ positive primitive representations.
\end {proof}

\begin {proof}[Proof of Theorem~\ref {th:etatwoaplusb}]
Let $z = \sigma (c)$, and $x = \sigma (a)$ and $y = \sigma (b)$ with the
purported generalized pentagonal numbers $a$ and $b$, where $\sigma$ is given
by \eqref{eq:sigmadef}. Then $c = 2 a + b$
translates as
\begin {equation}
\label {eq:eight}
z^2 + 2 = 2 x^2 + y^2,
\end {equation}
so we need to show that for $z \geq 11$ and coprime to~$6$, the integer
$k = z^2 + 2$ admits a positive representation $(x, y)$ by the quadratic
form $2 X^2 + Y^2$ other than $(x, y) = (1, z)$ and with $x$ and $y$
coprime to~$6$.

The existence of the primitive representation $(1, z)$ shows, using
Lemma~\ref {lm:eight} and the fact that $y$ is coprime to $6$, that all
prime divisors of~$k$ are $1$ or $3 \pmod 8$, and that as soon as $k$ has
at least two prime factors, there is another positive primitive
representation.
Notice that $k$ is divisible by~$3$, so we conclude that unless $k$ is a
power of~$3$, it admits a positive primitive representation $(x, y)$
with $x, y < z$. The following Proposition~\ref {prop:powersofthree} shows
that $k$ cannot be a power of~$3$ unless
$k = 3$ (and $z = 1$ and $c = 0$) or
$k = 27$ (and $z = 5$ and $c = 1$),
which are not covered by the theorem.

It remains to show that $x$ and $y$ can be taken coprime to~$6$.
Considering~\eqref {eq:eight} modulo~$8$ shows that $x$ and $y$ are
automatically odd.
The left hand side of~\eqref {eq:eight} is divisible by~$3$, while the
right hand side is divisible by~$3$ only if both $x$ and $y$ are coprime
to~$3$, or both are divisible by~$3$. The second possibility is ruled out
by the primitivity of the representation.
\end {proof}

\begin {proposition}
\label {prop:powersofthree}
The only solutions to $-2 = x^2 - 3^n$ with integers
$x$, $n \geq 0$ are given by $n=1$ and $x=1$, and by $n=3$ and $x=5$.
\end {proposition}

\begin {proof}
Assume that there are other solutions $(x, n)$ apart from the given ones.
If $n = 2m$ were even, then we would have
$\{ x - 3^m, x + 3^m \} \subseteq \{ \pm 1, \pm 2 \}$, whose only solution
is $x = 0$, $m = 0$. But this does not lead to a solution of the equation.
Write $n = 2 m + 1$ and let $y = 3^m$, so that
\begin {equation}
\label {eq:normsqrtthree}
-2 = x^2 - 3 y^2.
\end {equation}
Let $K = \Q (\sqrt 3)$. Then~\eqref {eq:normsqrtthree} is equivalent to
$x + y \sqrt 3$ being an element of $\OK = \Z [\sqrt 3]$ of norm~$-2$.
An initial solution is given by $\alpha = 1 + \sqrt 3$; according to
PARI/GP \cite {parigp} a fundamental unit of~$\OK$ is
$\epsilon = 2 + \sqrt 3$ of norm~$1$, so that all elements of $\OK$ of
norm~$-2$ are given by the $\pm \alpha \epsilon^k$ with $k \in \Z$.
Let $\rho : \sqrt 3 \mapsto - \sqrt 3$ denote the non-trivial Galois
automorphism of $K / \Q$. Since elements that are conjugate under~$\rho$
lead to the same solution of~\eqref {eq:normsqrtthree} up to the sign of~$y$,
$a^\rho = - a \epsilon^{-1}$ and
$(a \epsilon^{-k})^\rho = - a \epsilon^{k-1}$, it is enough to consider
solutions with $k \geq 0$ (which are in fact exactly the solutions with
$x$, $y > 0$).
Write $\alpha \epsilon^k = x_k + y_k \sqrt 3$ with $x_k$, $y_k \in \Z$.
Then
\[
x_0 = y_0 = 1, \quad
x_k = 2 x_{k-1} + 3 y_{k-1}
\text { and }
y_k = x_{k-1} + 2 y_{k-1}
\text { for } k \geq 1.
\]
To exclude the already known solutions with $n \leq 3$, we now switch to the
norm equation
\begin {equation}
\label {eq:normconductornine}
-2 = x^2 - 243 y^2
\end {equation}
in the order $\Oc = \Z [9 \sqrt 3]$ of conductor~$9$.
An initial solution is given by
$\alpha' = \alpha \epsilon^4 = 265 + 17 \cdot 9 \sqrt 3$,
and the fundamental unit of~$\Oc$ is
$\epsilon' = \epsilon^9 = 70226 + 4505 \cdot 9 \sqrt 3$, the smallest
power of~$\epsilon$ that lies in~$\Oc$.
Then the solutions of~\eqref {eq:normconductornine} (up to the signs of
$x$ and $y$) are derived from the
$\alpha' (\epsilon')^k = x_k + y_k \cdot 9 \sqrt 3$ with
$x_0 = 265$, $y_0 = 17$,
\[
x_k = 70226 x_{k-1} + 1094715 y_{k-1},
\quad
y_k = 4505 x_{k-1} + 70226 y_{k-1}
\text { for } k \geq 1.
\]
One notices that all $y_k$ are divisible by~$17$ and thus not a power of~$3$.
\end {proof}

\subsection {One multiplication}

The previous section gave an upper bound of one square and one multiplication
for each term of the series of~$\eta$. Even more favorable situations are
more difficult to analyze. They do not happen for all generalized pentagonal
numbers, and the non-existence of a primitive representation does not rule
out the existence of an imprimitive representation, which is enough for our
purposes and thus needs to be examined. For instance, the cases of one square
$c = 2 a$ or of one multiplication $c = a + b$ translate by
Lemma~\ref {lm:penta} into $z^2 + 1 = 2 x^2$ and $z^2 + 1 = x^2 + y^2$,
respectively, where $z = \sigma (c)$, $x = \sigma (a)$ and $y = \sigma (b)$.
Now $k = 2 x^2$ is the ``maximally imprimitive'' representation of
$k = x^2 + y^2$.

\begin {lemma}
\label {lm:four}
A positive integer~$k$ is primitively represented by the quadratic form
$X^2 + Y^2$ if and only if all its odd prime divisors are~$1 \pmod 4$ and
it is not divisible by~$4$. Its number of ordered positive primitive
representations is then given by $2^{\omega' (k) - 1}$.
\end {lemma}

The first part of the result is proved by Euler \cite{Euler58,Euler60}
using elementary arguments.
Concerning the number of representations, Euler \cite{Euler58} does not
provide a closed formula, but a number of arguments:
the case that $k$ is an odd prime is covered in~\S40;
factors of~$2$ are handled in~\S4;
the general case of odd square-free numbers follows by induction from~\S5,
where
``productum ex duobus huiusmodi numeris \textit {duplici modo} in duo quadrata
resolvi posse''
refers to the factor of~$2$ for each additional prime number, and examples
for products of two or three odd primes are given;
odd prime powers are not handled explicitly, but it should be possible to
derive the number of not necessarily primitive representations by induction
from the previous argument and then derive the number of primitive
representations from an inclusion--exclusion principle.
Again, modern number theory provides an easy proof.

\begin {proof}
The arguments are the same as in the proof of Lemma~\ref {lm:eight}, but
using the maximal order $\OK = \Z [i]$ of $K = \Q (i)$.
There are now four units $\{ \pm 1, \pm i \}$
instead of just two, but the unit $i$ only swaps $|x|$ and $|y|$, which is
taken into account by considering only ordered representations.
\end {proof}

\begin {theorem}
\label {th:etaaplusb}
A generalized pentagonal number $c \geq 2$ is the sum of two smaller ones,
that is, there are generalized pentagonal numbers $a, b < c$ such that
$c = a + b$, if and only if $12 c + 1$ is not a prime.
\end {theorem}

\begin {proof}
Let $z = \sigma (c)$, $x = \sigma (a)$ and $y = \sigma (b)$.
By Lemma~\ref {lm:penta}, $c = a + b$ is equivalent with $k = x^2 + y^2$
for $k = z^2 + 1 = 2 (12 c + 1)$, which is even, but not divisible by~$4$.
The existence of the primitive representation $(1, z)$ shows by
Lemma~\ref {lm:four} that all primes dividing $k / 2$ are $1 \pmod 4$,
and the lemma also implies that there is another primitive representation
unless $k = 2 p^\alpha$ with $p$~prime and $\alpha \geq 1$.
If $\alpha \geq 2$, we may take a primitive representation of $k / p^2$
and multiply it by~$p$. For $\alpha = 1$, there is no other representation.
\end {proof}

The first generalized pentagonal number that is not a sum of two previous
ones is $5 = 2 \cdot 2 + 1$. For larger numbers, it will be less and less
likely that $12 c + 1$ is prime. Heuristically, it is expected to happen
for only $O \left( \frac {\sqrt T}{\log T} \right)$ of the $\Theta (\sqrt T)$
generalized pentagonal numbers up to~$T$.

The first generalized pentagonal number requiring an imprimitive
representation is $c = 70$ with $z = 41$. From $41^2 + 1 = 2 \cdot 29^2$
we deduce $c = 2 a$ with the generalized pentagonal number $a = 35$.

Theorem~\ref {th:etaaplusb} proves the second point of
Theorem~\ref {th:addseq} for~$\eta$, since the $12 c + 1$ for generalized
pentagonal numbers $c$ are exactly the values of the two polynomials given
there, separately for ordinary pentagonal numbers and the other ones, and
omitting the single value $c = 0$.

\subsection {One squaring}

As seen in the previous section, it is possible that a generalized pentagonal
number is twice a previous one. But the following discussion shows that this
happens for a negligible (exponentially small) proportion of numbers.

By Lemma~\ref {lm:penta}, $c = 2 a$ translates into
$z^2 + 1 = 2 x^2$ for $z = \sigma (c)$ and $x = \sigma (a)$;
in other words, $z + x \sqrt 2$ is a unit of norm~$-1$ in
$\OK = \Z [\sqrt 2]$. An initial solution is given by the fundamental unit
$\epsilon = 1 + \sqrt 2$, of which exactly the odd powers
\[
\epsilon^{2 k + 1} = (1 + \sqrt 2)(3 + 2 \sqrt 2)^k = z_k + x_k \sqrt 2
\]
have norm~$-1$. They satisfy the linear recurrence
\[
z_0 = x_0 = 1, \quad
z_k = 3 z_{k-1} + 4 x_{k-1}, \quad
x_k = 2 z_{k-1} + 3 x_{k-1},
\]
which, considered modulo~$2$ and~$3$, shows that all the $z_k$ and $x_k$
are coprime to~$6$. However, growing exponentially, they are very rare.

\subsection {One cube}

In the cases where a generalized pentagonal number is not the sum of two
previous ones, it may still be three times a previous one, which leads to a
slightly faster computation of the term than by a square and a multiplication
according to Table~\ref {tab:cost}. But again, this case is exceedingly rare,
since $c = 3 a$ corresponds by Lemma~\ref {lm:penta} to
$z^2 + 2 = 3 x^2$ with $z = \sigma (c)$ and $x = \sigma (a)$.
Using the initial solution $z_0 = x_0 = 1$ and the fundamental unit
$2 + \sqrt 3$ of $\Z [\sqrt 3]$, all solutions are given by
\[
z_k = 2 z_{k-1} + 3 x_{k-1}, \quad x_k = z_{k-1} + 2 x_{k-1}.
\]
All $x_k$ and $z_k$ are odd, and $z_k \equiv (-1)^k \pmod 3$. However,
$3 \mid x_k$ for $k = 1$, or $k \geq 4$ and $4 \mid k$, in which cases
the associated~$a$ is not a generalized pentagonal number.

\section {Addition sequences for \texorpdfstring {$\theta$}{theta}-functions}
\label {sec:theta}

We now consider the Jacobi theta functions, showing that the
associated exponent sequences can be treated in analogy
with the pentagonal numbers for the eta function.

\subsection {Trigonal numbers and \texorpdfstring {$\theta_2$}{theta2}}
\label {ssec:thetatwo}

According to~\eqref {eq:thetatri}, the series for~$\theta_2$ can be computed
by an addition sequence for the \emph {trigonal numbers} $n (n + 1)$ for
$n \in \Z_{\geq 0}$. (The usual terminology calls the numbers $n(n + 1)/2$
\emph {triangular numbers} and excludes $n = 0$; the addition sequences for
triangular numbers are in bijection with those for trigonal numbers by
doubling each term of a sum and adding the initial step $2 = 1 + 1$.)

Trigonal numbers permit a characterization similar to that of generalized
pentagonal numbers in Lemma~\ref {lm:penta}:
The strictly increasing map $\sigma : c \mapsto \sqrt {4 c + 1}$ is a
bijection between trigonal numbers and odd positive integers.
So considering trigonal numbers $c$, $a$ and $b$ with
$z = \sigma (c)$, $x = \sigma (a)$ and $y = \sigma (b)$, we can write
$c = a + b$ if and only if $z^2 + 1 = x^2 + y^2$ and
$c = 2 a + b$ if and only if $z^2 + 2 = 2 x^2 + y^2$.
As for $\eta$, it is clear that there is an addition sequence for the trigonal
numbers with two additions per number using
\begin {align*}
a_0 & = 0 & a_n & = a_{n-1} + 2 = 2 n \\
b_0 & = 0 & b_n & = b_{n-1} + a_n = n (n + 1)
\end {align*}
The following result, which is analogous to Theorems~\ref {th:etaaplusb}
and~\ref {th:etatwoaplusb}, holds for trigonal numbers.

\begin {theorem}
\label {th:thetatwo}
A trigonal number $c \geq 6$ is the sum of two smaller ones if and only if
$2 c + 1$ is not a prime.
It is the sum of a smaller one and twice a smaller one if and only if
$4 c + 3$ is not a prime.
\end {theorem}

\begin {proof}
This follows from Lemma~\ref {lm:four} and~\ref {lm:eight}, using the same
techniques as in the proofs of Theorems~\ref {th:etaaplusb}
and~\ref {th:etatwoaplusb}. A subtlety arises for $c = 2 a + b$ when
$k = z^2 + 2 = p^\alpha = 2 x^2 + y^2$ is the power of a prime.
As there is no restriction on the divisibility of $z$ by~$3$, we may now
have $p \neq 3$. If $\alpha \geq 3$, the primitive representation for
$k / p^2$ can be multiplied by~$p$ as in the proof of
Theorem~\ref {th:etaaplusb}. If $\alpha = 2$, the primitive
representation $1 = 2 \cdot 0^2 + 1^2$ is degenerate and meaningless in
our context; then there is no second positive representation apart from
$k = 2 \cdot 1^2 + z^2$. However, $z^2 + 2 = p^2$ has no
solution in integers, so this case does in fact not occur.

As $z$ is odd, there is no such problem for $c = a + b$,
$k = z^2 + 1 = x^2 + y^2$, since then $k$ equals twice an odd number, and
even when $k = 2 p^2$ we can lift the primitive and positive representation
$2 = 1^2 + 1^2$.
\end {proof}

The addition sequence derived from Theorem~\ref{th:thetatwo} by letting $c = a + b$
whenever possible and $c = 2 a + b$ otherwise still has holes; the first
trigonal number~$c$ such that both $2 c + 1$ and $4 c + 3$ are prime is
$20 = 4 \cdot 5$. To fill these holes, one cannot use the generic addition
sequence above, as the sequence of the $a_n = 2 n$ is not contained in our
more optimized one. However, $20 = 12 + 6 + 2$, and the following general
result holds.

\begin {theorem}
\label {th:thetatwothree}
Every trigonal number $c \geq 6$ is the sum of at most three smaller ones.
\end {theorem}

\begin {proof}
Legendre has shown that every number is the sum of three triangular numbers
including~$0$ \cite [pp.\ 205--399]{Legendre97}. But this result is
useless in our context, since we do not wish to write a trigonal number
as a sum of itself and~$0$.
We need to solve $z^2 + 2 = x^2 + y^2 + t^2$ with odd $x$, $y$ and $t$.
The parity condition holds automatically from the fact that $z$ is odd,
as can be seen by examining the equation modulo~$4$. If only one of $x$,
$y$ and $t$ equals~$1$, we have found a meaningful representation of
$z^2 + 1 = x^2 + y^2$ and written the trigonal number as a sum of two
smaller ones. So we only need to show that there is another representation
of $z^2 + 2$ as a sum of three squares apart from the $24$ representations
obtained from $(x, y, t) = (z, 1, 1)$ by permutations or adding signs.
The number of primitive representations has been counted by
Gauss \cite [\S291]{Gauss01} \cite [Theorem~4.2]{Grosswald85},
for $k \geq 5$ and $k \equiv 3 \pmod 8$, as
$24 \, h (-k)$, where $h (-k)$ is the class number of the order of
discriminant~$-k$ in $\Q (\sqrt {-k})$. So we have an essentially different
primitive representation whenever $h (-k) \geq 2$, which is the case for
$c = (k^2 - 1) / 4 \geq 12$. For $c = 6$ we have $6 = 3 \cdot 2$,
corresponding to an imprimitive representation.
\end {proof}

Together, Theorems~\ref {th:thetatwo} and~\ref {th:thetatwothree} prove
the second point of Theorem~\ref {th:addseq} for $\theta_2$,
since the $2 c + 1$ for trigonal numbers are exactly the values of the
polynomial given there.

\subsection {Squares and \texorpdfstring {$\theta_0$}{theta0}
and \texorpdfstring {$\theta_1$}{theta1}}
\label {ssec:thetazero}

At first sight,
for the squares occurring as exponents of the usual series for $\theta_0$
and $\theta_1$, the relative scarcity of Pythagorean triples leaves little
hope of finding good addition sequences. Indeed, precise criteria are given
by Lemma~\ref {lm:four} and~\ref {lm:eight}. But whereas in
\S\S\ref {sec:eta} and~\ref {ssec:thetatwo} the existence of one primitive
representation was obvious from the shape of the numbers and we merely needed
to check whether a second, non-trivial representation existed,
in the case of squares there will be no primitive representation at all when
the number is divisible by a prime not satisfying the necessary congruences
modulo~$4$ or~$8$. However, Dobkin and Liption~\cite{DoLi80} show the existence of an addition
sequence for the first $N$ squares containing $N + O (N / \sqrt {\log N})$
terms by considering imprimitive representations; they also mention an
unpublished result, communicated by Donald Newman to Nicholas Pippenger,
that improves the bound to
$N + O (N / e^{c \log N / \log \log N})$ for some unknown constant $c > 0$.

Using a simple trick and the techniques of the previous sections, we may
easily obtain an asymptotically worse, but practically very satisfying
result, namely the second point of Theorem~\ref {th:addseq} for $\theta_0$
and $\theta_1$. For that, we split off one common
factor of~$q$ and consider exponents of the form $c = n^2 - 1$ for
$n \in \Z_{\geq 1}$, which we will call \emph {almost-square} in the
following.
The map $\sigma : c \mapsto \sqrt {c + 1}$ is a bijection between
almost-squares and positive integers.

\begin {theorem}
\label {th:thetazero}
An almost-square $c \geq 3$ is the sum of two smaller ones if and
only if $c + 2$ is neither a prime nor twice a prime.
It is the sum of a smaller one and twice a smaller one if and only if
$c + 3$ is neither a prime nor twice a prime nor twice the square of a prime.
\end {theorem}

\begin {proof}
The same techniques as for Theorem~\ref {th:thetatwo} apply.
As now we have no restriction any more on the parity of $z$ in
$k = z^2 + 1$ or $k = z^2 + 2$, we need to consider all the special cases
$k = p$, $k = 2 p$, $k = p^2$ (which cannot occur) and
$k = 2 p^2$ (which poses problems only for $k = 2 x^2 + y^2$ and
not for $k = x^2 + y^2$).
\end {proof}

\begin {theorem}
\label {th:thetazerothree}
Every almost-square $c \geq 24$ is the sum of at most three smaller
almost-squares.
\end {theorem}

\begin {proof}
The case $c$ even or equivalently $z = \sigma (c) = \sqrt {c + 1}$ odd
is handled as in Theorem~\ref {th:thetatwothree}, and we find a non-trivial
primitive representation for $z \geq 7$ with $c \geq 48$, and the
imprimitive representation $z^2 + 2 = 27 = 3 \cdot 3^2$ for $z = 5$ and
$c = 24 = 3 \cdot 8$.
In the case $c$ odd, $z$ even, the number of primitive representations is
given by Gauss as $12 \, h (- 4 k)$ for $k = z^2 + 2 = c + 3$, and we
have $h (- 4 k) \geq 3$ for $z \geq 6$.
\end {proof}

Together, Theorems~\ref {th:thetazero} and~\ref {th:thetazerothree} prove
the second point of Theorem~\ref {th:addseq} for $\theta_0$ and $\theta_1$.
If $c$ is an almost-square, then $c + 2$  can only be prime if
$c = (2 n)^2 - 1$ is odd, leading to the first polynomial of
Theorem~\ref {th:addseq}.
Conversely, $c + 2$ can only be twice a prime if $c = (2 n + 1)^2 - 1$
is even, leading to the second polynomial.

\subsection {Computing \texorpdfstring {$\theta$}{theta} functions
simultaneously}

The two previous sections have shown that good addition sequences for single
$\theta$ functions exist, which asymptotically approach an average of
one multiplication per term of the series (under the heuristic assumption
that the values of quadratic polynomials occurring in the theorems are
prime, or twice a prime, or twice the square of a prime with the same
logarithmic probabilities as arbitrary numbers).
In practice, one will often want to compute all $\theta$ functions
simultaneously. By considering all exponents at the same time, one may
potentially save a few additional multiplications.

Instead of considering almost-square numbers for $\theta_0$ and $\theta_1$,
we will revert to squares and consider the sequence \seqnum{A002620} of
\emph {quarter-squares}
$0$, $1$, $2$, $4$, $6$, $9$, $12$, $16$, $\ldots$
defined by
$t(n) = \lfloor (n+1)^2 / 4 \rfloor$ for $n \in \Z_{\geq 0}$,
which interleaves the squares $t(2m-1) = m^2$ and the
trigonal numbers $t(2m) = m(m+1)$ in increasing order.

\begin {theorem}
\label {th:thetatwoaplusb}
Every quarter-square $c > 1$ is the sum of a smaller one and twice a smaller one.
\end {theorem}

\begin {proof}
We use the following formula as a starting point:
\[
t (2 a n + \alpha)
= a^2 n^2 + a (\alpha + 1) n
+ \left\lfloor \frac {(\alpha + 1)^2}{4} \right\rfloor.
\]
Considering the primitive representation
$3^2 = 2 \cdot 2^2 + 1^2$,
it becomes natural to examine
\begin{multline}
t (6 n + \alpha) - 2 t (4 n + \beta) - t (2 n + \gamma)
= \\
(3 \alpha - 4 \beta - \gamma - 2) n
+ \left(
    \left\lfloor \frac {(\alpha + 1)^2}{4} \right\rfloor
- 2 \left\lfloor \frac {(\beta + 1)^2}{4} \right\rfloor
-   \left\lfloor \frac {(\gamma + 1)^2}{4} \right\rfloor
\right) .
\end{multline}
Then for each $\alpha \in \{ 0, \ldots, 5 \}$ there are $\beta$ and $\gamma$
for which this expression vanishes, and we obtain the following explicit
recursive formul{\ae} for the addition sequence:
\begin{align*}
t(6n)   &= 2 t(4n)   + t (2n-2) \\
t(6n+1) &= 2 t(4n)   + t (2n+1) \\
t(6n+2) &= 2 t(4n+1) + t (2n) \\
t(6n+3) &= 2 t(4n+2) + t (2n-1) \\
t(6n+4) &= 2 t(4n+2) + t (2n+2) \\
t(6n+5) &= 2 t(4n+3) + t (2n+1).
\end{align*}
\end {proof}

This shows that when computing all $\theta$ functions simultaneously, each
additional term of the series may be obtained with at most one squaring and
one multiplication, which has the merit of giving a uniform result without
any exceptions, but which is unfortunately worse than computing the
functions separately as in \S\S\ref {ssec:thetatwo} and~\ref {ssec:thetazero}
with only one multiplication per term most of the time.

To solve this problem, we consider yet another sequence of exponents
given by $f (n) = 2 \, \left\lfloor \frac {n^2}{8} \right\rfloor$
for $n \geq 1$, which interleaves in increasing order
the trigonal numbers $m (m+1)$ for $n = 2m+1$;
the even squares $(2 m)^2$ for $n = 4 m$;
and the even almost-squares, $(2 m + 1)^2 - 1$, for $n = 4 m + 2$.
Ignoring initial zeros, this sequence is equivalent to \seqnum{A182568}.
By separating the terms with odd and even exponents into two series and
by splitting off one power of~$q$ in the series with odd exponents, the
squares and almost-squares can be used to compute $\theta_0$ and $\theta_1$.

\begin {theorem}
\label {th:thetaaplusb}
Every element $c \geq 4$ in the sequence
$(f (n))_{n \geq 1} =
\left( 2 \left\lfloor \frac {n^2}{8} \right\rfloor \right)_{n \geq 1}$
is the sum of two smaller ones.
\end {theorem}

\begin {proof}
We may consider the sequence $g (n) = f (n) / 2$ in place of $f (n)$
itself. The starting point of the proof, which is similar to that of
Theorem~\ref {th:thetatwoaplusb}, is the following formula:
\[
g (4 a n + \alpha)
= 2 a^2 n^2 + a \alpha n + \left\lfloor \frac {\alpha^2}{8} \right\rfloor.
\]
We now replace $a$ by the elements of the Pythagorean triple
$5^2 = 4^2 + 3^2$ and compute
{
\small
\[
g (20 n + \alpha) - g (16 n + \beta) - g (12 n + \gamma)
= (5 \alpha - 4 \beta - 3 \gamma) n
+ \left(  \left\lfloor \frac {\alpha^2}{8} \right\rfloor
        - \left\lfloor \frac {\beta^2}{8} \right\rfloor
        - \left\lfloor \frac {\gamma^2}{8} \right\rfloor \right).
\]
}
It is easy to check that for every $\alpha \in \{ -9, \ldots, 10 \}$,
the values $\beta$ and $\gamma$ given in the following table make
this expression vanish.

\centerline {
\begin {tabular}{c|c|c}
$\alpha$ & $\beta$ & $\gamma$ \\
\hline
    $0$ &     $0$ &     $0$ \\
$\pm 1$ & $\pm 2$ & $\mp 1$ \\
$\pm 2$ & $\pm 1$ & $\pm 2$ \\
$\pm 3$ & $\pm 3$ & $\pm 1$ \\
$\pm 4$ & $\pm 2$ & $\pm 4$ \\
$\pm 5$ & $\pm 4$ & $\pm 3$ \\
$\pm 6$ & $\pm 6$ & $\pm 2$ \\
$\pm 7$ & $\pm 5$ & $\pm 5$ \\
$\pm 8$ & $\pm 7$ & $\pm 4$ \\
$\pm 9$ & $\pm 6$ & $\pm 7$ \\
   $10$ &     $8$ &     $6$
\end {tabular}
}
For $n = 0$ and $\alpha \in \{ 4, 6 \}$, the table entries lead to the
trivial relation $g (\alpha) = g (\alpha) + g (2)$, but one readily
verifies that $g (4) = 2 g (3)$ and $g (6) = 2 g (4)$.
\end {proof}

So when one or both of $\theta_0$ and $\theta_1$ are computed together with
$\theta_2$, the series may be evaluated with one multiplication per required
term, which proves the third point of Theorem~\ref {th:addseq}.

\section {Baby-step giant-step algorithm}
\label {sec:bsgs}

For evaluating the series expansions of the eta function and theta constants,
we may ignore the cost of multiplying by the
coefficients $c_n$ since they are all $1$ or $-1$.

To evaluate a power series truncated to include exponents $e_n \leq T$,
the baby-step giant-step algorithm of~\eqref{eq:bsgsformula} with
splitting parameter~$m$ requires
\begin{equation}
\label {eq:bsgscost}
(m-1) + (\lceil (T+1) / m \rceil - 1) \approx m + T / m
\end{equation}
multiplications.
The first term accounts for computing the powers
$q^2, \ldots q^m$ (baby-steps)
and the second term accounts for the multiplications by $q^m$ (giant-steps).
Setting $m \approx \sqrt{T}$ in \eqref{eq:bsgscost} gives the
minimized cost of $2 \sqrt{T} + O (1)$ multiplications.

The exponent sequences for the Jacobi theta functions
and the Dedekind eta function are just sparse enough so that
the baby-step giant-step algorithm performs worse
than computing the powers of~$q$ by an optimized addition sequence,
provided the latter is of length $N + o(N)$.
Indeed, there are $N \in \sqrt{T} + O (1)$ squares up to $T$, and
$N \in \sqrt{8T/3} + O (1) \approx 1.633 \sqrt{T} + O (1)$
generalized pentagonal numbers.

When computing all three theta functions simultaneously,
the baby-steps
can be recycled, but the giant-steps have to be done separately
for each function. The approximate cost of
\begin{equation}
\label {eq:bsgscost2}
m + 3 T / m
\end{equation}
is minimized by taking $m \in \sqrt{3 T} + O (1)$,
yielding $3.464 \sqrt{T} + O (1)$ multiplications.
This is again worse than computing the powers by an addition sequence,
since there are $2 \sqrt{T} + O (1)$ squares and trigonal numbers up to~$T$.
One gets slightly smaller constants for the baby-step giant-step
algorithm by recognizing that half of the powers $q, q^2, \ldots q^m$
can be computed using squarings, but the conclusion remains the same.

We can, however, do better in the baby-step giant-step algorithm
by choosing $m$ such that only a sparse subset of the exponents
$q^2, \ldots, q^m$ need to be computed.
For example, when considering squares $e_n = n^2$,
we seek $m$ such that there are few
squares modulo $m$. If we denote this number by $s(m)$, the
cost to minimize is
\begin{equation}
\label {eq:bsgscost3}
s(m)^{1+\varepsilon} + T / m.
\end{equation}
where the left term denotes the length of an addition sequence
for all the distinct values of $n^2 \bmod m$
as obtained, for instance, by Algorithm~\ref {alg:addseq}.
In the following, we show that $m$ can be chosen
so that \eqref{eq:bsgscost3} becomes $o(\sqrt{T})$,
giving an asymptotic speed-up.
In fact, Theorem~\ref{th:quadpolym} establishes this result
not only for squares, but for all quadratic exponent sequences.
We shall also explicitly derive suitable choices of $m$ for
squares, trigonal numbers, and generalized pentagonal numbers.

\subsection {Modular values of quadratic polynomials}
\label {ssec:modval}

\subsubsection{Squares}

We are interested in the number $s (m)$ of squares modulo a positive
integer~$m \geq 2$.
By the Chinese remainder theorem, $s (m)$ is a multiplicative number
theoretic function, so it is enough to consider the case
that $m = p^e$ is a power of some prime~$p$.
It is well-known that $(\Z / p^e \Z)^\times$ is cyclic of order
$p^{e-1} (p-1)$ if $p$ is odd, as shown by Gauss
\cite[\S\S52--56]{Gauss01} for $e=1$
and also by Gauss
\cite[\S\S82--89]{Gauss01} for $e \geq 2$;
that it is cyclic of order~$2^{e-1}$ if $p = 2$
and $e \in \{ 1, 2 \}$, and isomorphic to
$\Z / 2 \Z \times \Z / 2^{e-2}$ if $p = 2$ and $e \geq 3$,
again shown by Gauss \cite[\S\S90--91]{Gauss01}.
This determines the size of the kernel of the group endomorphism of
$(\Z / p^e \Z)^\times$ given by $x \mapsto x^2$, and shows that the size
of the image, that is, the number of squares modulo $p^e$ that are not
divisible by~$p$, is given by
$\frac {1}{2} p^{e-1} (p-1)$ if $p$ is odd;
by $1$ if $p = 2$ and $e \leq 3$;
and by $2^{e-3}$ if $p = 2$ and $e \geq 3$.

It remains to count the number of squares modulo~$p^e$ that are divisible
by~$p$. These are given by~$0$ and by the $p^{2k} z$, where $2 \leq 2k < e$
and $z$ is a square modulo~$p^{e - 2k}$ that is coprime to~$p$. So the number
of such squares is given by
\[
1 + \sum_{k = 1}^{\left\lfloor \frac {e-1}{2} \right\rfloor}
\left| \big( (\Z / p^{2k} \Z)^\times \big)^2 \right|.
\]
Distinguishing the cases that $p$ is odd or even, that $e$ is odd or even,
and using the result of the previous paragraph, a little computation
gives the total number of squares modulo $p^e$ as
\begin {equation}
\label {eq:squares}
s (p^e) =
\begin {cases}
\frac {1}{2} p^e - \frac {1}{2} p^{e-1}
+ \frac {p^{e-1} - p^{(e+1) \bmod 2}}{2 (p+1)} + 1,
& \text {for $p$ odd}; \\
2,
& \text {for $p = 2$ and $e \leq 2$}; \\
2^{e-3}
+ \frac {2^{e-3} - 2^{(e+1) \bmod 2}}{3} + 2,
& \text {for $p = 2$ and $e \geq 3$};
\end {cases}
\end {equation}
where the exponent $(e+1) \bmod 2$ is understood to be~$0$ or~$1$.

We are interested in low numbers of squares, that is,
small values of the ratio $s (m) / m$.
Let $p_k$ denote the $k$-th prime and let $\theta (x)$ denote the logarithm
of the product of all primes not exceeding~$x$.
Then~\eqref {eq:squares} shows that the sequence $s (m) / m$
tends to~$0$ roughly as $1 / 2^k$ for $m = e^{\theta (p_k)}$,
so that the inferior limit of the full sequence of $s (m) / m$ is~$0$.
We consider the subsequence of ratios providing successive minima,
in the sense that $s (m) / m < s (m') / m'$ for all $m' < m$;
the $m$ realizing these successive minima are given by the sequence
\seqnum{A085635},
the corresponding $s (m)$ form sequence \seqnum{A084848}.
Using~\eqref {eq:squares} and the multiplicativity of~$s (m)$,
one readily computes the values of these sequences for $m \leq 10^8$,
see Table~\ref {tab:squaresminima}; we have augmented the table by
the values for the $m = e^{\theta (p_k)}$.

\begin {table}
\begin {center}
\begin {footnotesize}
\begin {tabular}{rr@{\,}c@{\,}lrl}
$k$ & \multicolumn {3}{l}{$m = \mathrm {\seqnum{A085635}} (k)$}
    & $s (m) = \mathrm {\seqnum{A084848}} (k)$
    & $s (m) / m$ \\
\hline
$  1$ & $       2$ & $=$
 & $2$
 & $     2$ & $1.0   $ \\
$  2$ & $       3$ & $=$
 & $3$
 & $     2$ & $0.67  $ \\
$  3$ & $       4$ & $=$
 & $2^2$
 & $     2$ & $0.50  $ \\
      & $       6$ & $=$
 & $2 \cdot 3$
 & $     4$ & $0.67  $ \\
$  4$ & $       8$ & $=$
 & $2^3$
 & $     3$ & $0.38  $ \\
$  5$ & $      12$ & $=$
 & $2^2 \cdot 3$
 & $     4$ & $0.33  $ \\
$  6$ & $      16$ & $=$
 & $2^4$
 & $     4$ & $0.25  $ \\
      & $      30$ & $=$
 & $2 \cdot 3 \cdot 5$
 & $    12$ & $0.40  $ \\
$  7$ & $      32$ & $=$
 & $2^5$
 & $     7$ & $0.22  $ \\
$  8$ & $      48$ & $=$
 & $2^4 \cdot 3$
 & $     8$ & $0.17  $ \\
$  9$ & $      80$ & $=$
 & $2^4 \cdot 5$
 & $    12$ & $0.15  $ \\
$ 10$ & $      96$ & $=$
 & $2^5 \cdot 3$
 & $    14$ & $0.15  $ \\
$ 11$ & $     112$ & $=$
 & $2^4 \cdot 7$
 & $    16$ & $0.14  $ \\
$ 12$ & $     144$ & $=$
 & $2^4 \cdot 3^2$
 & $    16$ & $0.11  $ \\
      & $     210$ & $=$
 & $2 \cdot 3 \cdot 5 \cdot 7$
 & $    48$ & $0.23  $ \\
$ 13$ & $     240$ & $=$
 & $2^4 \cdot 3 \cdot 5$
 & $    24$ & $0.10  $ \\
$ 14$ & $     288$ & $=$
 & $2^5 \cdot 3^2$
 & $    28$ & $0.097 $ \\
$ 15$ & $     336$ & $=$
 & $2^4 \cdot 3 \cdot 7$
 & $    32$ & $0.095 $ \\
$ 16$ & $     480$ & $=$
 & $2^5 \cdot 3 \cdot 5$
 & $    42$ & $0.088 $ \\
$ 17$ & $     560$ & $=$
 & $2^4 \cdot 5 \cdot 7$
 & $    48$ & $0.086 $ \\
$ 18$ & $     576$ & $=$
 & $2^6 \cdot 3^2$
 & $    48$ & $0.083 $ \\
$ 19$ & $     720$ & $=$
 & $2^4 \cdot 3^2 \cdot 5$
 & $    48$ & $0.067 $ \\
$ 20$ & $    1008$ & $=$
 & $2^4 \cdot 3^2 \cdot 7$
 & $    64$ & $0.063 $ \\
$ 21$ & $    1440$ & $=$
 & $2^5 \cdot 3^2 \cdot 5$
 & $    84$ & $0.058 $ \\
$ 22$ & $    1680$ & $=$
 & $2^4 \cdot 3 \cdot 5 \cdot 7$
 & $    96$ & $0.057 $ \\
$ 23$ & $    2016$ & $=$
 & $2^5 \cdot 3^2 \cdot 7$
 & $   112$ & $0.056 $ \\
      & $    2310$ & $=$
 & $2 \cdot 3 \cdot 5 \cdot 7 \cdot 11$
 & $   288$ & $0.12  $ \\
$ 24$ & $    2640$ & $=$
 & $2^4 \cdot 3 \cdot 5 \cdot 11$
 & $   144$ & $0.055 $ \\
$ 25$ & $    2880$ & $=$
 & $2^6 \cdot 3^2 \cdot 5$
 & $   144$ & $0.050 $ \\
$ 26$ & $    3600$ & $=$
 & $2^4 \cdot 3^2 \cdot 5^2$
 & $   176$ & $0.049 $ \\
$ 27$ & $    4032$ & $=$
 & $2^6 \cdot 3^2 \cdot 7$
 & $   192$ & $0.048 $ \\
&&& $\vdots$ \\
$ 94$ & $41801760$ & $=$
 & $2^5 \cdot 3^2 \cdot 5 \cdot 7 \cdot 11 \cdot 13 \cdot 29$
 & $211680$ & $0.0051$ \\
$ 95$ & $42325920$ & $=$
 & $2^5 \cdot 3^2 \cdot 5 \cdot 7 \cdot 13 \cdot 17 \cdot 19$
 & $211680$ & $0.0050$ \\
$ 96$ & $48454560$ & $=$
 & $2^5 \cdot 3^2 \cdot 5 \cdot 7 \cdot 11 \cdot 19 \cdot 23$
 & $241920$ & $0.0050$ \\
$ 97$ & $49008960$ & $=$
 & $2^6 \cdot 3^2 \cdot 5 \cdot 7 \cdot 11 \cdot 13 \cdot 17$
 & $217728$ & $0.0044$ \\
$ 98$ & $54774720$ & $=$
 & $2^6 \cdot 3^2 \cdot 5 \cdot 7 \cdot 11 \cdot 13 \cdot 19$
 & $241920$ & $0.0044$ \\
$ 99$ & $61261200$ & $=$
 & $2^4 \cdot 3^2 \cdot 5^2 \cdot 7 \cdot 11 \cdot 13 \cdot 17$
 & $266112$ & $0.0043$ \\
$100$ & $68468400$ & $=$
 & $2^4 \cdot 3^2 \cdot 5^2 \cdot 7 \cdot 11 \cdot 13 \cdot 19$
 & $295680$ & $0.0043$ \\
$101$ & $82882800$ & $=$
 & $2^4 \cdot 3^2 \cdot 5^2 \cdot 7 \cdot 11 \cdot 13 \cdot 23$
 & $354816$ & $0.0043$ \\
$102$ & $89535600$ & $=$
 & $2^4 \cdot 3^2 \cdot 5^2 \cdot 7 \cdot 11 \cdot 17 \cdot 19$
 & $380160$ & $0.0042$
\end {tabular}
\end {footnotesize}
\end {center}
\caption {\label {tab:squaresminima}
Successive minima of $s (m) / m$ for squares.}
\end {table}

\subsubsection {Trigonal numbers}
We now turn to general quadratic polynomials
$a X^2 + b X + c \in \Z [X]$.
Completing the square as
$a \left( X + \frac {b}{2 a} \right)^2
+ \left( c - \frac {b^2}{4 a} \right)$
shows that they take as many values modulo~$m$ as there are squares,
unless $\gcd (m, 2 a) \neq 1$; and hereby, rational coefficients
$a$, $b$, $c$ are also permitted as long as the denominators are coprime
to~$m$.

For the polynomial $X^2 + X$ defining trigonal numbers, this means that
their number of values modulo odd prime powers is still given
by~\eqref {eq:squares}.
Modulo $2^e$, one quickly verifies that $x$ and $a - x$ yield the
same value of $X^2 + X$ if and only if
$(a + 1)(2 x - a) \equiv 0 \pmod {2^e}$, which yields the exact two
solutions $a = -1$ (for $a$ odd) and $a = 2 x$ (for $a$ even).
So $X^2 + X$ takes each value twice or zero times; as all its values
are even, it takes the even values exactly twice, and there are $2^{e - 1}$
of them.

Table~\ref {tab:trigonalsminima} summarizes the successive minima of the
ratio between the number $t (m)$ of trigonal numbers modulo~$m$ and $m$
for $m \leq 10^4$ and some values just below $10^8$.

\begin {table}
\begin {center}
\begin {footnotesize}
\begin {tabular}{rr@{\,}c@{\,}lrl}
$k$ & \multicolumn {3}{l}{$m$}
    & $t (m)$
    & $t (m) / m$ \\
\hline
$  1$ & $       2$ & $=$
 & $2$
 & $     1$ & $0.50  $ \\
$  2$ & $       6$ & $=$
 & $2 \cdot 3$
 & $     2$ & $0.33  $ \\
$  3$ & $      10$ & $=$
 & $2 \cdot 5$
 & $     3$ & $0.30  $ \\
$  4$ & $      14$ & $=$
 & $2 \cdot 7$
 & $     4$ & $0.29  $ \\
$  5$ & $      18$ & $=$
 & $2 \cdot 3^2$
 & $     4$ & $0.22  $ \\
$  6$ & $      30$ & $=$
 & $2 \cdot 3 \cdot 5$
 & $     6$ & $0.20  $ \\
$  7$ & $      42$ & $=$
 & $2 \cdot 3 \cdot 7$
 & $     8$ & $0.19  $ \\
$  8$ & $      66$ & $=$
 & $2 \cdot 3 \cdot 11$
 & $    12$ & $0.18  $ \\
$  9$ & $      70$ & $=$
 & $2 \cdot 5 \cdot 7$
 & $    12$ & $0.17  $ \\
$ 10$ & $      90$ & $=$
 & $2 \cdot 3^2 \cdot 5$
 & $    12$ & $0.13  $ \\
$ 11$ & $     126$ & $=$
 & $2 \cdot 3^2 \cdot 7$
 & $    16$ & $0.13  $ \\
$ 12$ & $     198$ & $=$
 & $2 \cdot 3^2 \cdot 11$
 & $    24$ & $0.12  $ \\
$ 13$ & $     210$ & $=$
 & $2 \cdot 3 \cdot 5 \cdot 7$
 & $    24$ & $0.11  $ \\
$ 14$ & $     330$ & $=$
 & $2 \cdot 3 \cdot 5 \cdot 11$
 & $    36$ & $0.11  $ \\
$ 15$ & $     390$ & $=$
 & $2 \cdot 3 \cdot 5 \cdot 13$
 & $    42$ & $0.11  $ \\
$ 16$ & $     450$ & $=$
 & $2 \cdot 3^2 \cdot 5^2$
 & $    44$ & $0.098 $ \\
$ 17$ & $     630$ & $=$
 & $2 \cdot 3^2 \cdot 5 \cdot 7$
 & $    48$ & $0.076 $ \\
$ 18$ & $     990$ & $=$
 & $2 \cdot 3^2 \cdot 5 \cdot 11$
 & $    72$ & $0.073 $ \\
$ 19$ & $    1170$ & $=$
 & $2 \cdot 3^2 \cdot 5 \cdot 13$
 & $    84$ & $0.072 $ \\
$ 20$ & $    1386$ & $=$
 & $2 \cdot 3^2 \cdot 7 \cdot 11$
 & $    96$ & $0.069 $ \\
$ 21$ & $    1638$ & $=$
 & $2 \cdot 3^2 \cdot 7 \cdot 13$
 & $   112$ & $0.068 $ \\
$ 22$ & $    2142$ & $=$
 & $2 \cdot 3^2 \cdot 7 \cdot 17$
 & $   144$ & $0.067 $ \\
$ 23$ & $    2310$ & $=$
 & $2 \cdot 3 \cdot 5 \cdot 7 \cdot 11$
 & $   144$ & $0.062 $ \\
$ 24$ & $    2730$ & $=$
 & $2 \cdot 3 \cdot 5 \cdot 7 \cdot 13$
 & $   168$ & $0.062 $ \\
$ 25$ & $    3150$ & $=$
 & $2 \cdot 3^2 \cdot 5^2 \cdot 7$
 & $   176$ & $0.056 $ \\
$ 26$ & $    4950$ & $=$
 & $2 \cdot 3^2 \cdot 5^2 \cdot 11$
 & $   264$ & $0.053 $ \\
$ 27$ & $    5850$ & $=$
 & $2 \cdot 3^2 \cdot 5^2 \cdot 13$
 & $   308$ & $0.053 $ \\
$ 28$ & $    6930$ & $=$
 & $2 \cdot 3^2 \cdot 5 \cdot 7 \cdot 11$
 & $   288$ & $0.042 $ \\
$ 29$ & $    8190$ & $=$
 & $2 \cdot 3^2 \cdot 5 \cdot 7 \cdot 13$
 & $   336$ & $0.041 $ \\
&&& $\vdots$ \\
$107$ & $47477430$ & $=$
 & $2 \cdot 3^2 \cdot 5 \cdot 7 \cdot 11 \cdot 13 \cdot 17 \cdot 31$
 & $290304$ & $0.0061$ \\
$108$ & $49639590$ & $=$
 & $2 \cdot 3^2 \cdot 5 \cdot 7 \cdot 11 \cdot 13 \cdot 19 \cdot 29$
 & $302400$ & $0.0061$ \\
$109$ & $51482970$ & $=$
 & $2 \cdot 3^2 \cdot 5 \cdot 7 \cdot 11 \cdot 17 \cdot 19 \cdot 23$
 & $311040$ & $0.0060$ \\
$110$ & $60090030$ & $=$
 & $2 \cdot 3^2 \cdot 5 \cdot 7 \cdot 11 \cdot 13 \cdot 23 \cdot 29$
 & $362880$ & $0.0060$ \\
$111$ & $60843510$ & $=$
 & $2 \cdot 3^2 \cdot 5 \cdot 7 \cdot 13 \cdot 17 \cdot 19 \cdot 23$
 & $362880$ & $0.0060$ \\
$112$ & $76715730$ & $=$
 & $2 \cdot 3^2 \cdot 5 \cdot 7 \cdot 13 \cdot 17 \cdot 19 \cdot 29$
 & $453600$ & $0.0059$ \\
$113$ & $82006470$ & $=$
 & $2 \cdot 3^2 \cdot 5 \cdot 7 \cdot 13 \cdot 17 \cdot 19 \cdot 31$
 & $483840$ & $0.0059$ \\
$114$ & $87297210$ & $=$
 & $2 \cdot 3^3 \cdot 5 \cdot 7 \cdot 11 \cdot 13 \cdot 17 \cdot 19$
 & $498960$ & $0.0057$ \\
$115$ & $95611230$ & $=$
 & $2 \cdot 3^2 \cdot 5 \cdot 11 \cdot 13 \cdot 17 \cdot 19 \cdot 23$
 & $544320$ & $0.0057$
\end {tabular}
\end {footnotesize}
\end {center}
\caption {\label {tab:trigonalsminima}
Successive minima of $t (m) / m$ for trigonal numbers.}
\end {table}

\subsubsection {Generalized pentagonal numbers}
The number of values $p (m)$ of the polynomial $\frac {(3 X - 1) X}{2}$
modulo~$m$ is given by~\eqref {eq:squares} outside of~$2$ and~$3$.

Modulo~$3^e$, it is a bijection. If it takes the same value at
$x$ and $x + a$, then $a (6 x + 3 a - 1) \equiv 0 \pmod {3^e}$,
so $a = 0$.

Modulo powers of~$2$, it is to be understood that the values of the
polynomial in~$\Q [X]$ in integer arguments, which are integers, are reduced
modulo~$2^e$. The number of such values equals the number of values of
$(3 X - 1) X$ modulo $2^{e+1}$. As with the trigonal numbers examined above,
this polynomial takes every even value twice: It takes the same value in $x$
and in $a - x$ if and only if
$(3 a - 1)(a - 2 x) \equiv 0 \pmod {2^{e+1}}$,
which has the even solution $a = 2 x$ and the odd solution
$a = 1/3 \bmod {2^{e+1}}$.
So the number of values is $2^e$, and the polynomial induces a bijection
of $\Z / 2^e \Z$.

Successive minima of $p (m) / m$ for $m$ up to $10^4$ and just below $10^8$
are given in Table~\ref {tab:pentagonalsminima}. In line with the previous
reasoning, none of the $m$ achieving a successive minimum is divisible
by~$2$ or~$3$.

\begin {table}
\begin {center}
\begin {footnotesize}
\begin {tabular}{rr@{\,}c@{\,}lrl}
$k$ & \multicolumn {3}{l}{$m$}
    & $p (m)$
    & $p (m) / m$ \\
\hline
$  1$ & $       2$ & $=$
 & $2$
 & $     2$ & $1.0   $ \\
$  2$ & $       5$ & $=$
 & $5$
 & $     3$ & $0.60  $ \\
$  3$ & $       7$ & $=$
 & $7$
 & $     4$ & $0.57  $ \\
$  4$ & $      11$ & $=$
 & $11$
 & $     6$ & $0.55  $ \\
$  5$ & $      13$ & $=$
 & $13$
 & $     7$ & $0.54  $ \\
$  6$ & $      17$ & $=$
 & $17$
 & $     9$ & $0.53  $ \\
$  7$ & $      19$ & $=$
 & $19$
 & $    10$ & $0.53  $ \\
$  8$ & $      23$ & $=$
 & $23$
 & $    12$ & $0.52  $ \\
$  9$ & $      25$ & $=$
 & $5^2$
 & $    11$ & $0.44  $ \\
$ 10$ & $      35$ & $=$
 & $5 \cdot 7$
 & $    12$ & $0.34  $ \\
$ 11$ & $      55$ & $=$
 & $5 \cdot 11$
 & $    18$ & $0.33  $ \\
$ 12$ & $      65$ & $=$
 & $5 \cdot 13$
 & $    21$ & $0.32  $ \\
$ 13$ & $      77$ & $=$
 & $7 \cdot 11$
 & $    24$ & $0.31  $ \\
$ 14$ & $      91$ & $=$
 & $7 \cdot 13$
 & $    28$ & $0.31  $ \\
$ 15$ & $     119$ & $=$
 & $7 \cdot 17$
 & $    36$ & $0.30  $ \\
$ 16$ & $     133$ & $=$
 & $7 \cdot 19$
 & $    40$ & $0.30  $ \\
$ 17$ & $     143$ & $=$
 & $11 \cdot 13$
 & $    42$ & $0.29  $ \\
$ 18$ & $     175$ & $=$
 & $5^2 \cdot 7$
 & $    44$ & $0.25  $ \\
$ 19$ & $     275$ & $=$
 & $5^2 \cdot 11$
 & $    66$ & $0.24  $ \\
$ 20$ & $     325$ & $=$
 & $5^2 \cdot 13$
 & $    77$ & $0.24  $ \\
$ 21$ & $     385$ & $=$
 & $5 \cdot 7 \cdot 11$
 & $    72$ & $0.19  $ \\
$ 22$ & $     455$ & $=$
 & $5 \cdot 7 \cdot 13$
 & $    84$ & $0.18  $ \\
$ 23$ & $     595$ & $=$
 & $5 \cdot 7 \cdot 17$
 & $   108$ & $0.18  $ \\
$ 24$ & $     665$ & $=$
 & $5 \cdot 7 \cdot 19$
 & $   120$ & $0.18  $ \\
$ 25$ & $     715$ & $=$
 & $5 \cdot 11 \cdot 13$
 & $   126$ & $0.18  $ \\
$ 26$ & $     935$ & $=$
 & $5 \cdot 11 \cdot 17$
 & $   162$ & $0.17  $ \\
$ 27$ & $    1001$ & $=$
 & $7 \cdot 11 \cdot 13$
 & $   168$ & $0.17  $ \\
$ 28$ & $    1309$ & $=$
 & $7 \cdot 11 \cdot 17$
 & $   216$ & $0.17  $ \\
$ 29$ & $    1463$ & $=$
 & $7 \cdot 11 \cdot 19$
 & $   240$ & $0.16  $ \\
$ 30$ & $    1547$ & $=$
 & $7 \cdot 13 \cdot 17$
 & $   252$ & $0.16  $ \\
$ 31$ & $    1729$ & $=$
 & $7 \cdot 13 \cdot 19$
 & $   280$ & $0.16  $ \\
$ 32$ & $    1925$ & $=$
 & $5^2 \cdot 7 \cdot 11$
 & $   264$ & $0.14  $ \\
$ 33$ & $    2275$ & $=$
 & $5^2 \cdot 7 \cdot 13$
 & $   308$ & $0.14  $ \\
$ 34$ & $    2975$ & $=$
 & $5^2 \cdot 7 \cdot 17$
 & $   396$ & $0.13  $ \\
$ 35$ & $    3325$ & $=$
 & $5^2 \cdot 7 \cdot 19$
 & $   440$ & $0.13  $ \\
&&& $\vdots$ \\
$128$ & $76491415$ & $=$
 & $5 \cdot 7 \cdot 11 \cdot 13 \cdot 17 \cdot 29 \cdot 31$
 & $1088640$ & $0.014 $ \\
$129$ & $80925845$ & $=$
 & $5 \cdot 7 \cdot 11 \cdot 13 \cdot 19 \cdot 23 \cdot 37$
 & $1149120$ & $0.014 $ \\ $130$ & $82944785$ & $=$
 & $5 \cdot 7 \cdot 11 \cdot 17 \cdot 19 \cdot 23 \cdot 29$
 & $1166400$ & $0.014 $ \\
$131$ & $88665115$ & $=$
 & $5 \cdot 7 \cdot 11 \cdot 17 \cdot 19 \cdot 23 \cdot 31$
 & $1244160$ & $0.014 $ \\
$132$ & $98025655$ & $=$
 & $5 \cdot 7 \cdot 13 \cdot 17 \cdot 19 \cdot 23 \cdot 29$
 & $1360800$ & $0.014 $
\end {tabular}
\end {footnotesize}
\end {center}
\caption {\label {tab:pentagonalsminima}
Successive minima of $p (m) / m$ for pentagonal numbers.}
\end {table}

In the remainder of this section, we let $f (m)$ denote the number of distinct
values modulo~$m$ taken by a quadratic polynomial~$F$ (including, but not
limited to, the number of squares, trigonal and generalized pentagonal numbers
modulo~$m$). Our goal is to prove that a judicious choice of~$m$ leads to
$f (m) / m$ sufficiently small so that the baby-step giant-step algorithm
applied to a $q$-series with exponents given by~$F$ (and trivial coefficients)
takes sublinear time in the number of terms of the series. We use standard
estimates of analytic number theory, as given, for instance, in \cite{RoSc62},
and elementary analytic arguments like the following observation.

\begin {lemma}
\label {lm:analytic}
Consider functions $k$, $m : \N \to \N$. If
\[
| \log m (n) - k (n) \log k (n) | \in o (k (n) \log k (n))
\text { and }
k (n) \to \infty \text { for } n \to \infty,
\]
then
\[
k \in \Theta (\log m / \log \log m).
\]
\end {lemma}

\begin {proof}
More precisely, we show that
$k (n) \log \log m (n) / \log m (n) \to 1$ for $n \to \infty$,
so the constant implied in $\Theta$-notation is~$1$.
The main hypothesis can be reformulated as
\begin {equation}
\label {eq:limit}
\frac {\log m}{k \log k} \to 1.
\end {equation}
Since $\frac {x}{y} \to 1$ implies $\frac {\log x}{\log y} \to 1$ for
$y$ positive and bounded away from~$0$, we also have
\[
\frac {\log \log m}{\log k \left( 1 + \frac {\log \log k}{\log k} \right)}
\to 1.
\]
With $k \to \infty$, we obtain $\frac {\log \log m}{\log k} \to 1$,
and the desired statement follows from a division by~\eqref {eq:limit}.
\end {proof}

\begin {theorem}
\label {th:quadpolym}
For a fixed quadratic polynomial $F (X) \in \Q [X]$ that takes
integral values at integral arguments, and for
$N \to \infty$, a judicious choice of~$m$ (detailed in the proof)
leads to an effective constant $c > 0$
such that the baby-step giant-step algorithm computes the series
$\sum_{n=1}^N q^{F (n)}$ with
$N^{1 - c / \log \log N}$
multiplications, which grows more slowly than
$N / \log^r N$ for any $r > 0$.
\end {theorem}

\begin {proof}
The number of multiplications of the baby-step giant-step algorithm is
bounded above by
\begin {equation}
\label {eq:babygiantsteps}
\frac {F (N)}{m} + 2 f (m) \log_2 m + O (1)
\in
O \left( \frac {N^2}{m} + f (m) \log m \right),
\end {equation}
where the first term accounts for the giant-steps and
$f (m)$ is the number of values of $F$ modulo~$m$.
Here we pessimistically assume that each of the values is obtained
by a separate addition chain using at most $\log_2 m$ doublings and
$\log_2 m$ additions; in practice, we expect the number of additional
terms in the addition sequence to be negligible and the number of
multiplications to be rather of the order of~$f (m)$.

Following the discussion for trigonal numbers above, we have $f (m) = s (m)$
if $m$ is odd and coprime to the common denominator of the coefficients of~$F$
and to its leading coefficient. To minimize~$s (m)$ and thus~$f (m)$,
following~\eqref {eq:squares} it becomes desirable to build~$m$ with as many
prime factors as possible.
So we choose~$m$ as the product of the first $k$ primes $p_1, \ldots, p_k$,
but avoiding~$2$ and the finitely many primes dividing the leading coefficient
and the denominator of~$F$.
For the time being, $k$ is an unknown function of the desired number
of terms~$N$; it will be fixed later to minimize~\eqref {eq:babygiantsteps}.
The quantity $m$ depends on~$k$ (and thus ultimately also on~$N$) and on~$F$.
We will have $k (N) \to \infty$ and $m \to \infty$ as $N \to \infty$.

In a first step, we estimate $k$ in terms of $m$, uniformly for all~$F$.
Let $\theta (x)$ denote the logarithm
of the product of all primes not exceeding~$x$.
For $N \to \infty$, the finite number of primes excluded from~$m$ have a
negligible impact, so that
\begin {equation}
\label {eq:mpk}
| \log m - \theta (p_k) | \in O (1).
\end {equation}
We use the following standard results \cite[Theorems~4~and~3]{RoSc62} from analytic number theory:
\begin {equation}
\label {eq:rstheta}
| \theta (x) - x| \in O (x / \log x)
\end {equation}
and
\begin {equation}
\label {eq:pk}
| p_k - k \log k | \in O (k \log \log k).
\end {equation}

 From \eqref {eq:pk} we obtain
$\frac {p_k}{k \log k} \to 1$ and, as in the proof of Lemma~\ref {lm:analytic},
$\frac {\log p_k}{\log k} \to 1$, so that
$\frac {p_k}{\log p_k} \in \Theta (k)$ after division.
Together with~\eqref {eq:rstheta}, in which $x$ has been replaced by $p_k$,
this implies
\begin {equation}
\label {eq:thetapk}
| \theta (p_k) - p_k | \in O (k).
\end {equation}
Summing up \eqref {eq:mpk}, \eqref {eq:thetapk} and~\eqref {eq:pk}, and using the
triangle inequality implies
\[
| \log m - k \log k |
\in O (k \log \log k).
\]
Lemma~\ref {lm:analytic} then implies that $k \in \Theta (\log m / \log \log m)$,
so that $p_k \in \Theta (\log m)$ by~\eqref {eq:pk}.
(In other words, the largest prime contributes roughly $\log \log m$
bits, so that $\log m / \log \log m$ primes are needed.)

Now by~\eqref {eq:squares} we have
\[
f (m) = s (m) = \prod_{p \mid m} \frac {p+1}{2}
\in
O \left(
\frac {m}{2^k}
\prod_{i=1}^k \left( 1 + \frac {1}{p_i} \right)
\right)
\subseteq
O \left(
\frac {m \log \log m}{2^k}
\right),
\]
where the last inclusion stems from
\[
0 \leq
\log \left( \prod_{i=1}^k \left( 1 + \frac {1}{p_i} \right) \right)
\leq \sum_{i=1}^k \frac {1}{p_i}
\in \log \log p_k + \Theta (1)
\]
by \cite[Theorem~5]{RoSc62},
and from the above relation between $p_k$ and $m$.

So the second term of~\eqref {eq:babygiantsteps} lies in
$O \left( m \log m \log \log m / 2^k \right)$.
We now use $k \in \Theta (\log m / \log \log m)$, and
let $c_1 > 0$ be such that
$2^k \geq m^{2 c_1 / \log \log m}$.
Since $\log m \log \log m \in O \left( m^{c_1 / \log \log m} \right)$,
the second term of~\eqref {eq:babygiantsteps} is an element of
$O \left( m^{1 - c_1 / \log \log m} \right)$.

We may still choose the magnitude of $m$ with respect to~$N$.
We let $m \in \Theta \left( N^{1 + c_2 / \log \log N} \right)$
for a sufficiently small $0 < c_2 \approx c_1 / 2$,
so that the first term of the complexity~\eqref {eq:babygiantsteps}
lies in $O \left( N^{1 - c_2 / \log \log N} \right)$.
Moreover, $| \log \log m - \log \log N | \in o (1)$, so the second term
is essentially bounded by
$N^{1 + (c_2 - c_1) / \log \log N} \approx N^{1 - c_2 / \log \log N}$;
in any case, there is a $0 < c_3 < c_2$ such that the second term lies
in $O \left( N^{1 - c_3 / \log \log N} \right)$.
By replacing $c_3$ with a suitable $0 < c < c_3$, the constant of the
big-Oh can be made~$1$ (or any other positive value).
\end {proof}

\subsection{Implementation}

To realize the baby-step giant-step algorithm for
computing a sum such as $\sum_{n^2 \leq T} q^{n^2}$, we may
use a precomputed table of word-sized $m$ for which $s(m)/m$
attains its successive minima, and a table of corresponding values $s(m)$.

Given $T$, we search the table to choose the $m$ minimizing $T / m + s(m)$.
Once $m$ is chosen, we create a table of the baby-step
exponents (the squares modulo~$m$)
and insert by Algorithm~\ref{alg:addseq} extra exponents into this table
as necessary until its entries form an addition sequence.
Few such insertions are needed in practice, making $s(m)$
an accurate estimate for the length of the baby-step addition sequence,
unlike the pessimistic bound $2 s(m) \log_2(m)$ in the proof
of Theorem~\ref{th:quadpolym}.

\begin{figure}[hbt]
    \centering
    \includegraphics[width=0.6\textwidth]{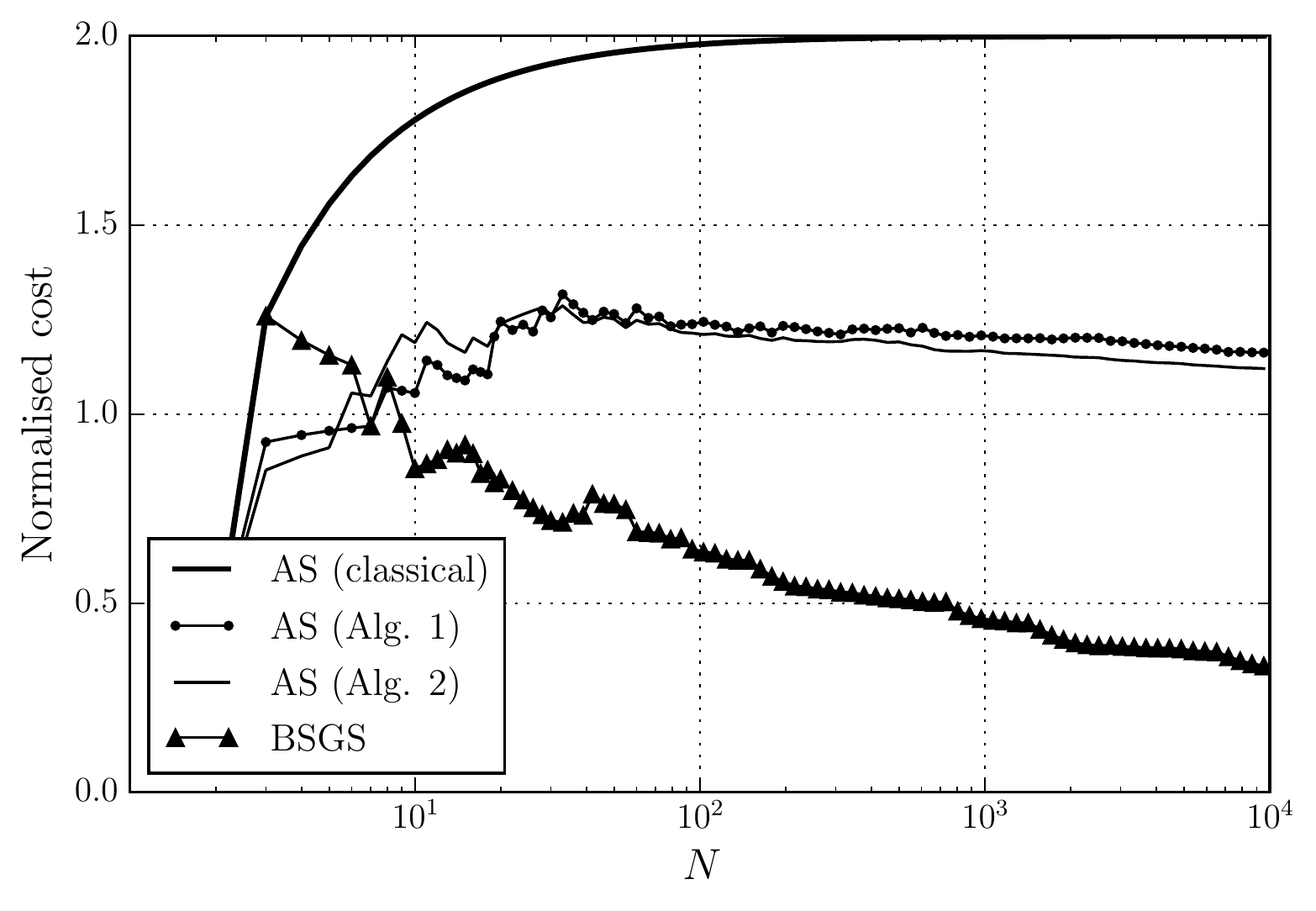}
    \caption{Normalized theoretical cost in the FFT model
           ($(3m + 2.333s) / (3N)$ for $m$
           complex multiplications and $s$ complex squarings)
           to add the first $N$ nonzero terms in the $q$-series of the
           Dedekind eta function, using three different addition sequences (AS)
           or the baby-step giant-step algorithm (BSGS).
           The classical addition sequences approaches~$2$ multiplications
           per term.
           The short addition sequences generated with
           Algorithm~\ref{alg:addseq} and Algorithm~\ref{alg:eta} both
           approach~$1$ multiplication per term.
           BSGS is asymptotically better than any addition sequence.}
    \label {fig:etaefficiency}
\end{figure}

The procedure is, of course, analogous with trigonal numbers
or generalized pentagonal numbers as exponents.
Figure~\ref{fig:etaefficiency} illustrates the theoretical
speed-up for generalized pentagonal numbers based on counting multiplications.

\section{Benchmarks}
\label {sec:bench}

We have implemented the Jacobi theta functions and the
Dedekind eta function in the Arb library~\cite{Arb},
and complemented the existing implementation of the Dedekind eta function
in release 0.3 of CM~\cite {cm} by the baby-step giant-step algorithm.
The complete function evaluation
involves three steps:
\begin{enumerate}
\item Reduce $\tau \in \C, \Im (\tau) > 0$
to the fundamental domain of the action of $\Sl_2 (\Z)$ on the upper complex
half-plane.
\item Compute $q = \exp(2 \pi i \tau)$ (for $\eta$) or
$q = \exp(\pi i \tau)$ (for $\theta_i$).
\item Sum the truncated $q$-series.
\end{enumerate}
The first step has negligible cost.
For the third step,
we have implemented the optimized addition sequences (AS)
as well as the baby-step giant-step (BSGS) algorithm.
AS and BSGS are used automatically at low and high precision, respectively.
Here, we compare both methods.
The measurements were done on an Intel Core i5-4300U CPU running
a 64-bit Linux kernel. MPIR~2.7.2~\cite{mpir} was used for multiprecision arithmetic.
Tables~\ref{tab:timeeta}, \ref{tab:timethetasimul} and \ref{tab:timethetasingle}
show timings with Arb,
and Table~\ref{tab:timecomparison} compares old implementations with the
new implementations in Arb and CM.

We take
$\tau = (-B + \sqrt{D}) / (2A)$ with $A = 1305, B = 1523, D = -6961631$,
which is a typical complex multiplication point occurring in class
polynomial construction.
The magnitude $|q| \approx 0.00174$ is slightly
smaller than the worst case $|q| \approx 0.00433$ at the corner of the
fundamental domain.
For $p$ bits of floating point precision, the truncation order is
$T \approx 0.11 p$ when computing
the eta function and $T \approx 0.22 p$ when computing theta functions.

Our code includes a small practical optimization:
For a tolerance of $2^{-p}$, the term $q^n$ needs to be computed
to a precision of only $p - n |\log_2(|q|)|$ bits,
and we change the internal precision for each term accordingly.
Empirically, this saves roughly a factor of~$1.5$ when using addition
sequences and a factor of~$1.2 $ in the BSGS algorithm (which
is improved less since the baby-steps have to be done at
essentially full precision). The speed-up of the BSGS algorithm
compared to addition sequences is therefore somewhat smaller
than what one might predict by counting multiplications.

\begin {table}[ht]
\footnotesize
\begin {center}
\begin {tabular}{r@{\;}rd{9}d{9}d{9}rr}
Bits    & $T$ & \multicolumn {1}{c}{Exponential} &
\multicolumn {1}{c}{Sum (AS)} & \multicolumn {1}{c}{Sum (BSGS)} &
Speed-up & Theory \\
\hline
$10^2$  &       7 & 0.000\,001\,59 & 0.000\,001\,79 & 0.000\,002\,86 & 0.63 & 0.74 \\
$10^3$  &     100 & 0.000\,018\,0  & 0.000\,026\,1  & 0.000\,023\,6  & 1.11 & 1.34 \\
$10^4$  &    1080 & 0.001\,52      & 0.001\,69      & 0.001\,20      & 1.41 & 1.63 \\
$10^5$  &   10880 & 0.066\,1       & 0.128          & 0.080\,9       & 1.58 & 2.06 \\
$10^6$  &  108676 & 1.74           & 6.12           & 3.11           & 1.97 & 2.32 \\
$10^7$  & 1090987 & 32.4           & 259            & 119            & 2.18 & 2.77 \\
\end {tabular}
\caption {\label {tab:timeeta}
Timings for computing the Dedekind eta function
at different precisions. From left to right: bit precision,
truncation order $T$ (last included exponent),
seconds to compute $q = \exp(2 \pi i \tau)$,
seconds to evaluate the sum using (AS),
seconds to evaluate the sum using (BSGS),
measured speed-up AS / BSGS,
and theoretical speed-up based on counting multiplications
in the FFT cost model.}
\end {center}
\end {table}

Timings for the eta function are shown in Table~\ref{tab:timeeta}.
Here, AS is the optimized addition sequence of Algorithm~\ref{alg:eta}.
We time the complex exponential separately and only show the speed-up
of (BSGS) over (AS) for the series summation.
At lower precision, the complex exponential takes a comparable
amount of time to summing the series, making the real-world speed-up smaller.
The speed-up for the series summation by itself is nonetheless useful
since there are situations where $q$ is available
without the need to compute the full complex exponential,
for example, during batch evaluation over an arithmetic progression
of $\tau$ values.

\begin {table}[ht]
\begin {center}
\begin {tabular}{rrd{9}d{9}rrr}
Bits    & $T$     &
\multicolumn {1}{c}{Sum (AS)} & \multicolumn {1}{c}{Sum (BSGS)} &
Speed-up & Theoretical \\
\hline
$10^2$  &      20 & 0.000\,003\,50 & 0.000\,005\,80 & 0.60 & 0.67 \\
$10^3$  &     210 & 0.000\,038\,6  & 0.000\,049\,2  & 0.78 & 0.89 \\
$10^4$  &    2162 & 0.002\,29      & 0.002\,16      & 1.06 & 1.18 \\
$10^5$  &   21756 & 0.178          & 0.134          & 1.33 & 1.55 \\
$10^6$  &  218089 & 8.97           & 5.71           & 1.57 & 1.78 \\
$10^7$  & 2181529 & 380            & 199            & 1.91 & 2.18 \\
\end {tabular}
\caption {\label {tab:timethetasimul} Time in seconds
to compute the theta functions $\theta_0(\tau), \theta_1(\tau), \theta_2(\tau)$
simultaneously, given $q = e^{\pi i \tau}$. Timings to compute the
complex exponential are the same as in Table~\ref{tab:timeeta} and thus omitted.}
\end {center}
\end {table}

\begin {table}[ht]
\begin {center}
\begin {tabular}{rrd{9}d{9}rr}
Bits    & $T$     &
\multicolumn {1}{c}{Sum (AS)} & \multicolumn {1}{c}{Sum (BSGS)} &
Speed-up & Theoretical \\
\hline
$10^2$  &      16 & 0.000\,002\,44 & 0.000\,003\,21 & 0.76 & 0.84 \\
$10^3$  &     196 & 0.000\,031\,1  & 0.000\,024\,5  & 1.27 & 1.51 \\
$10^4$  &    2116 & 0.001\,86      & 0.001\,10      & 1.69 & 2.23 \\
$10^5$  &   21609 & 0.147          & 0.065\,3       & 2.25 & 2.88 \\
$10^6$  &  218089 & 6.81           & 2.67           & 2.55 & 2.95 \\
$10^7$  & 2181529 & 280            & 90.1           & 3.11 & 3.58 \\
\end {tabular}
\caption {\label {tab:timethetasingle} Time in seconds
to compute $\theta_0(\tau)$ alone, given $q = e^{\pi i \tau}$.
Timings to compute the
complex exponential are the same as in Table~\ref{tab:timeeta} and thus omitted.}
\end {center}
\end {table}

Table~\ref{tab:timethetasimul} shows the corresponding timings to
compute three theta functions simultaneously.
Here, AS uses Theorem~\ref{th:thetaaplusb}, computing each term with
one multiplication or squaring.
The speed-up for BSGS is
smaller compared to the eta function
since three independent giant-step evaluations are done,
one for each theta function.
Table~\ref{tab:timethetasingle} shows timings for
computing a single theta function. For AS, we use Algorithm~\ref{alg:addseq}
to generate a short addition sequence for the squares alone.
Here the BSGS algorithm gives the largest speed-up.

Finally, Table~\ref{tab:timecomparison} compares the new implementations
with three previous implementations for the evaluation of $\eta(\tau)$.
The \texttt{eta} function in PARI/GP~\cite {parigp}
uses the classical addition sequence without the precision trick.
CM~\cite {cm} in version 0.2.1 used an optimized addition sequence
(Algorithm~\ref{alg:eta}) without the precision trick.
An implementation of the AGM method courtesy of R\'egis Dupont
(unpublished, cf.~\cite {Dupont11})
is also tested.
All implementations were linked
against the same version 2.7.2 of MPIR for multiprecision arithmetic.

\begin {table}[ht]
\begin {center}
\begin {tabular}{r@{\:}|d{6}d{6}d{6}|d{6}d{6}}
Bits    & \multicolumn {1}{c}{PARI/GP} & \multicolumn {1}{c}{CM-0.2.1} &
\multicolumn {1}{c|}{AGM} & \multicolumn {1}{c}{New CM-0.3} & \multicolumn {1}{c}{New Arb} \\
\hline
$10^4$  & 0.008\,69 & 0.004\,57 & 0.007\,89 & 0.002\,97 & 0.002\,72 \\
$10^5$  & 0.654     & 0.284     & 0.245     & 0.164     & 0.147   \\
$10^6$  & 29.9      & 10.9      & 6.78      & 4.77      & 4.85    \\
$10^7$  & 1\,310    & 440       & 124       & 150       & 151
\end {tabular}
\caption {\label {tab:timecomparison} Time in seconds to compute $\eta(\tau)$.}
\end {center}
\end {table}

In conclusion, we achieve a small but measurable speed-up. At practical
precisions, the baby-step giant-step algorithm saves somewhat less than a
factor of~$2$ in running time over an optimized addition sequence, which
itself saves a factor of~$2$ over the widely used classical addition sequence.
Using an optimized addition sequence instead of the classical addition
sequence raises the crossover point
between series evaluation and the AGM method from about $10^4$ bits to $10^5$ bits,
while the baby-step giant-step method raises it further to about $10^6$ bits,
very roughly.
The exact crossover point will vary depending on the system, the
libraries used for multiprecision arithmetic, the size of $|q|$
(a smaller value is more advantageous for series evaluation),
and whether $q$ needs to be computed from $\tau$ by a full exponential
evaluation.

\section {Acknowledgments}
This research was partially funded by ERC Starting Grant ANTICS 278537
and DFG Priority Project SPP 1489.

\bibliographystyle {plain}
\bibliography {addseq}

\bigskip
\hrule
\bigskip

\noindent 2010 {\it Mathematics Subject Classification}: Primary 11Y55; Secondary 11B83, 11F20.

\noindent \emph{Keywords:}
addition sequence, pentagonal number, theta function, Dedekind eta function, polynomial evaluation, baby-step giant-step algorithm

\bigskip
\hrule
\bigskip

\begin {flushleft}
\noindent (Concerned with sequences
\seqnum{A001318},
\seqnum{A002620},
\seqnum{A084848},
\seqnum{A085635}, and
\seqnum{A182568}.)
\end {flushleft}

\end {document}